\documentclass{ifacconf}

\usepackage{natbib}        

\usepackage{amsmath,amssymb,amsfonts}
\usepackage[ruled,vlined,algo2e]{algorithm2e}
\usepackage{graphicx}
\usepackage{textcomp}
\usepackage{stmaryrd}
\usepackage{tabularx}
\usepackage{multirow}
\usepackage{float}

\usepackage[font=footnotesize]{caption}
\usepackage{subcaption}

\usepackage{color}

\newcommand{\tnf}[1]{\textnormal{#1}}
\newcommand{\mbs}[1]{\boldsymbol{#1}}

\newcommand{\mbf}[1]{\mathbf{ #1}}

\newcommand{\mcl}[1]{\mathcal{#1}}

\newcommand{\R}{\mathbb{R}}
\newcommand{\N}{\mathbb{N}}

\newcommand{\ip}[2]{\left\langle #1, #2 \right\rangle}

\newcommand{\from}[0]{\leftarrow}

\usepackage{upgreek}
\newcommand{\PI}[0]{\Uppi}
\newcommand{\PIset}[0]{\mathbf{\Pi}}

\newcommand{\bmat}[1]{\begin{bmatrix}#1\end{bmatrix}}

\newcommand{\smallbmat}[1]{\left[\scriptsize\begin{smallmatrix}
		#1\end{smallmatrix} \right]}
\newcommand{\mat}[1]{\begin{matrix}#1\end{matrix}}

\newcommand{\lbmat}[1]{\left[\!\!\begin{array}{l}#1\end{array}\!\!\right]}

\newcommand{\slbmat}[1]{\small\left[\!\begin{array}{l}#1\end{array}\!\right]}

\newcommand{\slmat}[1]{\small\begin{array}{l}#1\end{array}}

\floatstyle{ruled}
\newfloat{block}{thbp}{lop}
\floatname{block}{Block}

\let\bl\bigl
\let\bbl\Bigl
\let\bbbl\biggl

\let\br\bigr
\let\bbr\Bigr
\let\bbbr\biggr

\begin{document}
	
\begin{frontmatter}
	
	\title{
		$H_{\infty}$-Optimal Estimator Synthesis for Coupled Linear 2D PDEs using Convex Optimization\thanksref{footnoteinfo}%
	}
	\thanks[footnoteinfo]{This work was supported by the NSF grant CMMI-1935453.}
	
	\author[First]{Declan S. Jagt} 
	\author[First]{Matthew M. Peet} 
	
	\address[First]{Arizona State University, (e-mail: djagt@asu.edu, mpeet@asu.edu).}
	
	\begin{abstract}
		Any suitably well-posed PDE in two spatial dimensions can be represented as a Partial Integral Equation (PIE) -- with system dynamics parameterized using Partial Integral (PI) operators. Furthermore, $L_2$-gain analysis of PDEs with a PIE representation can be posed as a linear operator inequality, which can be solved using convex optimization. In this paper, these results are used to derive a convex-optimization-based test for constructing an $H_{\infty}$-optimal estimator for 2D PDEs. In particular, we first use PIEs to represent an arbitrary well-posed 2D PDE where sensor measurements occur along some boundary of the domain. An associated Luenberger-type estimator is then parameterized using a PI operator $\mcl{L}$ as the observer gain. Examining the error dynamics of this estimator, it is proven that an upper bound on the $H_{\infty}$-norm of these error dynamics can be minimized by solving a linear operator inequality on PI operator variables. Finally, an analytical formula is proposed for inversion of a class of 2D PI operators, which is then used to reconstruct the Luenberger gain $\mcl{L}$. Results are implemented in the PIETOOLS software suite -- applying the methodology and simulating the resulting observer for an unstable 2D heat equation with boundary observations.
	\end{abstract}
	
	\begin{keyword}
		Distributed Parameter Systems, Observer Synthesis, PDEs, LMIs.
	\end{keyword}
	
\end{frontmatter}

	\maketitle
	\thispagestyle{empty}
	
	
	
	\vspace*{-0.2cm}
	\section{INTRODUCTION}	\vspace*{-0.2cm}

Partial Differential Equations (PDEs) are frequently used to model physical systems, relating the temporal evolution of an internal state variable to its spatial distribution. For example, to model the density $\mbf{u}(t,x,y)$ of an exponentially growing, distributed population on a spatial domain $(x,y)\in[0,1]^2$, we can use the following 2D PDE (see e.g.~\cite{holmes1994partial})
\begin{align}\label{eq:PDE_example_intro}
	\mbf{u}_{t}(t)&=\mbf{u}_{xx}(t)+\mbf{u}_{yy}(t)+r\mbf{u}(t) +w(t),	\\
	z(t)&=\int_{0}^{1}\int_{0}^{1}\mbf{u}(t,x,y)dx dy,	\notag
\end{align}
wherein $r$ is a parameter determining the population growth, $w(t)$ is an external disturbance, $z(t)$ denotes the total population size, and where the evolution of the state is further constrained by boundary conditions such as
\begin{equation*}
	\mbf{u}(t,0,y)=\mbf{u}_{x}(t,1,y)= 0,\qquad \mbf{u}(t,x,0)=\mbf{u}_{y}(t,x,1)= 0.
\end{equation*}
For state feedback control of such systems, we require real-time knowledge of the distributed internal state $\mbf{u}(t)$. However, in practice, direct measurement of the distributed state would require a prohibitive number of sensors. To alleviate the sensing burden, therefore, we commonly make a smaller number of observations -- typically on the boundary of the domain.
For example, in the population model, we might measure population density on the upper boundaries ($x=1$ and $y=1$), yielding observed outputs
\begin{equation*}
	\mbf{q}_{1}(t,y)=\mbf{u}(t,1,y)\qquad \text{and}\qquad \mbf{q}_{2}(t,x)=\mbf{u}(t,x,1).
\end{equation*}
The role of an estimator, then, is to take these limited observations and use them to reconstruct the distributed state at all points in its domain. Unfortunately, the infinite-dimensional nature of the system dynamics significantly complicates the problem of estimator design.

For comparison, consider a linear Ordinary Differential Equation (ODE), with a finite-dimensional state $u(t)\in\R^{n}$ and observed output $q(t)\in\R^{m}$,
\begin{align*}
	\dot{u}(t)&=Au(t)+Bw(t),	&
	q(t)&=C_{q}u(t)+D_{q}w(t).
\end{align*}
The most common approach to estimating the state, $u(t)$, is to construct a Luenberger-type observer, with state estimate $\hat{u}(t)$ and parameterized by a gain matrix $L$ as
\begin{equation*}
	\dot{\hat{u}}(t)=A\hat{u}(t) +Bw(t) +L(C_{q}\hat{u}(t)-q(t)).
\end{equation*}
Now, define a regulated output $z(t)$ and estimate $\hat{z}(t)$ of the regulated output as
\begin{equation*}
	z(t)=C u(t)+Dw(t),\qquad \hat{z}(t)=C\hat{u}(t).
\end{equation*}
Then the matrix $L$ which minimizes $\sup_{w\in L_2\setminus\{0\}}\frac{\|\hat{z}-z\|_{L_2}}{\|w\|_{L_2}}$ (the $H_{\infty}$-norm) may be found by solving the Linear Matrix Inequality (LMI)
\begin{align*}
	\min_{\gamma>0,P,W}\gamma,&	\\[-2.15em]
	&	&
	\!\!\bmat{-\gamma I&-D&C\\(\cdot)^*&-\gamma I&-PB^*-WD_{q}\\(\cdot)^*&(\cdot)^*&(\cdot)^*+PA+WC_{q}}&\!\preceq 0, \\[-1.25em]
	\text{s.t.}\enspace P\succ 0,&
\end{align*}
and setting $L=P^{-1}W$  (see e.g.~\cite{duan2013LMIs}). 

However, consider now a 2D PDE such as in~\eqref{eq:PDE_example_intro} with state $\mbf{u}(t,x,y)$.
Observing the value of this state along the boundary of the domain (e.g. $\mbf{q}_{1}(t,y)=\mbf{u}(t,1,y)$), the resulting sensed output will be infinite-dimensional (a 1D function). 
It is relatively simple to define an equivalent of the Luenberger-type estimator for ODEs, where we have
\begin{align*}
	\hat{\mbf{u}}_{t}(t)&=\hat{\mbf{u}}_{xx}(t)+\hat{\mbf{u}}_{yy}(t)+r\hat{\mbf{u}}(t)+w(t)+\mcl{L}(\hat{\mbf{q}}(t)-\mbf{q}(t)), \\
	\hat{\mbf{q}}_{1}(t,y)&=\hat{\mbf{u}}(t,1,y),\qquad \hat{\mbf{q}}_{2}(t,x)=\hat{\mbf{u}}(t,x,1),
\end{align*}
with regulated output estimate $\hat{z}(t)=\int_{0}^{1}\int_{0}^{1}\hat{\mbf{u}}(t,x,y)dx dy$.
The goal, then, is to find the observer gain, $\mcl{L}$, which minimizes the $H_{\infty}$-norm of the map from $w$ to $\hat z-z$. However this requires us to parameterize infinite-dimensional operators and optimize performance of PDE systems.

To avoid 
the problems of parameterizing infinite-dimensional operators and characterizing performance 
of PDE systems, a common approach is to project the PDE state onto a finite-dimensional subspace -- yielding a linear ODE -- and synthesizing an estimator based on this finite-dimensional approximation. 
Recent applications of this approach include: 1D systems with observer delay in~\cite{lhachemi2021outputfeedackstabilization}, 1D stochastic systems in~\cite{wang2024Observer_Stochastic}, and 2D Navier-Stokes equations in~\cite{zayats2021ObserverNS}, each deriving LMI conditions for verifying stability of the resulting error dynamics. 
However, parameterizing an estimator only by a finite-dimensional operator (a matrix) necessarily introduces conservatism.
In addition, properties such as optimality of the estimator for the ODE do not a priori guarantee optimality or even convergence of the estimator for the PDE, see e.g.~\cite{zuazua2005Observability_FiniteDifference}. As such, conditions for convergence of the constructed estimator must be proven a posteriori.


Aside from projection methods, perhaps the most common approach for estimator synthesis of PDEs is the backstepping method. 
Using this approach, a Luenberger-type estimator is parameterized by a multiplier operator, and convergence of the estimator is ensured by mapping the resulting error dynamics to a stable target system, using an integral transformation.
Using this approach, estimators can be designed for a variety of 1D PDEs, including e.g. semi-linear parabolic systems in~\cite{meurer2013BacksteppingObserver_Semilinear}, hyperbolic systems in~\cite{yu2020BacksteppingObserver_traffic}, and ODE-PDE cascade systems in~\cite{hasan2016BacksteppingObserver_PDEODE}, as well as PDEs in multiple spatial variables as in~\cite{jadachowski2015BacksteppingObserver_ND}. 
A disadvantage of this approach, however, is that the observer gains are defined by multiplier operators, thus introducing conservatism. Furthermore, each estimator is constructed only for a narrow class of systems -- and extending the approach to new systems may require significant expertise. In addition, the backstepping method offers no guarantee of optimality of the obtained estimators.


In this paper, we propose an alternative, LMI-based method for constructing a Luenberger-type estimator for a general class of linear, 2nd-order, 2D PDEs. In particular, we focus on systems of the form
\begin{align*}
	\mbf{u}_{t}(t)&=\textstyle{\sum_{i,j=0}^{2}}\:A_{ij}\partial_{x}^{i}\partial_{y}^{j}\mbf{u}(t) +Bw(t),\qquad \mbf{u}(t)\in X,	\\
	z(t)&=\textstyle{\sum_{i,j=0}^{2}}\int_{[0,1]^2}[C_{ij}]\partial_{x}^{i}\partial_{y}^{j}\mbf{u}(t) +Dw(t),	
\end{align*}
where $X\subseteq L_2^{n}[[0,1]^2]$ is the domain of the PDE, defined by a set of linear boundary conditions.
In addition, we assume that the value of the state is observed along the boundary of the domain, yielding an infinite-dimensional sensed output $\mbf{q}(t)\in L_2^{m}[0,1]$. For example, we allow for a sensed output along the upper boundary $y=1$ as
\begin{align*}
	\mbf{q}(t,x)=\textstyle{\sum_{i,j=0}^{1}}\:F_{ij}(x)\partial_{x}^{i}\partial_{y}^{j}\mbf{u}(t,x,1) +Gw(t).
\end{align*}
We allow for similar outputs on other boundaries as well.
To construct an estimator for the system with such outputs, we adopt an approach similar to that presented for 1D PDEs with finite-dimensional sensed outputs $\mbf{q}(t)\in\R^{m}$ in~\cite{das2019PIE_estimation}. In that paper, the 1D PDE was first converted to an equivalent Partial Integral Equation (PIE), expressing the dynamics of the system in terms of a \textit{fundamental state} $\mbf{v}(t)$ associated to the PDE as
\begin{align*}
	\mcl{T}\mbf{v}_{t}(t)&=\mcl{A}\mbf{v}(t)+\mcl{B}w(t),	&
	z(t)&=\mcl{C}\mbf{v}(t)+\mcl{D}w(t),	\\
	\mbf{q}(t)&=\mcl{C}_{q}\mbf{v}(t) +\mcl{D}_{q}w(t),
\end{align*}
where now $z(t)$ and $\mbf{q}(t)$ are still the regulated and sensed outputs, respectively, and where the parameters ($\mcl{T}$, $\mcl{A}$, etc.) are all Partial Integral (PI) operators. Given this PIE representation of the system, the authors then proposed constructing a Luenberger-type estimator as
\begin{align*}
	\mcl{T}\hat{\mbf{v}}_{t}(t)&=\mcl{A}\hat{\mbf{v}}(t) +\mcl{B}w(t) +\mcl{L}(\mcl{C}_{q}\hat{\mbf{v}} -\mbf{q}(t)),	&
	\hat{z}(t)&=\mcl{C}\hat{\mbf{v}}(t).
\end{align*}
Finally, the authors showed that if there exist PI operators $\mcl{P}$ and $\mcl{W}$ that solve the linear operator inequality
\begin{align}\label{eq:LPI_intro}
	\min_{\gamma>0,\mcl{P},\mcl{W}}&\gamma,	\notag\\[-2.15em]
	&	&
	&\bmat{-\gamma I &\! -\mcl{D} &\! \mcl{C}\\
		(\cdot)^* &\! -\gamma I &\! -[\mcl{B}^*\mcl{P}+\mcl{D}_{q}^*\mcl{W}^*]\mcl{T}	\\
		(\cdot)^* &\! (\cdot)^* &\! (\cdot)^* + \mcl{T}^*[\mcl{P}\mcl{A}+\mcl{W}\mcl{C}_{q}]}	
	\preceq 0,\notag\\[-1.25em]	
	\text{s.t.}\enspace \mcl{P}\succ& \: 0,
\end{align}
then, using the estimator gain $\mcl{L}=\mcl{P}^{-1}\mcl{W}$, the $H_{\infty}$-norm of the associated map from $w$ to $\hat{z}-z$ is upper-bounded by $\gamma$.
This linear operator inequality can be efficiently solved using convex optimization methods with the PIETOOLS software suite (\cite{shivakumar2021PIETOOLS}).

Unfortunately, when using the approach presented in~\cite{das2019PIE_estimation} to synthesize estimators for systems of 2D PDEs, we encounter several challenges. In particular, the sensed outputs from a 2D PDE are not finite-dimensional. This requires a more complicated parameterization of the gain $\mcl{L}:L_2^{m}[0,1]\to L_2^{n}[[0,1]^2]$ and hence the variable $\mcl{W}$ in~\eqref{eq:LPI_intro} in order to avoid conservatism. 
Furthermore, 
although it has been shown that a broad class of 2D PDEs with finite-dimensional inputs and outputs can be represented as 2D PIEs in~\cite{jagt2022PIE_2DHinfty}, 
a similar representation has not been derived for systems involving infinite-dimensional outputs. 
This poses the challenge of deriving a PIE representation for the sensed outputs $\mbf{q}(t)\in L_2^{m}[0,1]$, in terms of the fundamental state $\mbf{v}(t)\in L_2^{n}[[0,1]^2]$.
Finally, although a Luenberger-type estimator for the resulting 2D PIE representation may again be synthesized by solving the operator inequality in~\eqref{eq:LPI_intro}, computing the associated estimator gain $\mcl{L}=\mcl{P}^{-1}\mcl{W}$ requires inverting the PI operator $\mcl{P}$ -- raising the question of how to compute the inverse of PI operators in 2D. 

In the remainder of this paper, we address each of these challenges. 
First, in Subsec.~\ref{sec:main_result:PIE}, we derive an equivalent PIE representation for a broad class of 2D PDEs with infinite-dimensional sensed outputs. In Subsec.~\ref{sec:main_result:LPI}, we parameterize $\mcl{W}$ in~\eqref{eq:LPI_intro}, and show how optimal estimator synthesis for the PIE can be performed by solving the operator inequality. Finally, in Sec.~\ref{sec:LPI_implementation}, we derive an explicit expression for the inverse of a class of 2D PI operators, and show how the operator inequality in~\eqref{eq:LPI_intro} can be solved as an LMI. In Sec.~\ref{sec:Numerical_Examples}, we implement the methodology via the software suite PIETOOLS, and use numerical simulation to illustrate the approach for an unstable heat equation.

	\section{Preliminaries}
	
	\subsection{Notation}
	
	For a given domain $\Omega\subset\R^2$, let $L_2^n[\Omega]$ denote the set of $\R^n$-valued square-integrable functions on $\Omega$, where we omit the domain when clear from context. 
 	Define $W_{2}^n[[0,1]^2]$ as a Sobolev subspace of $L_2^n[[0,1]^2]$, where
 	\begin{equation*}
 		W_{2}^n[[0,1]^2]\!=\!\bl\{\mbf{v}~\bl|\: \partial_x^{\alpha_1}\partial_y^{\alpha_2}\mbf{v}\!\in\! L_2^n[[0,1]^2],\: \forall \alpha\!\in\!\N^2\!:\! \alpha_{1},\alpha_{2}\!\leq 2\br\}\!.
 	\end{equation*}
	For $\mbf{v}\in W_{2}^n[[0,1]^2]$, we denote the Dirac delta operators
	\begin{equation*}
		(\Delta_{x}^{0} \mbf{v})(y):=\mbf{v}(0,y)\quad \text{and}\quad (\Delta_{y}^{1} \mbf{v})(x):=\mbf{v}(x,1).
	\end{equation*}
	For a function $N\in L_2^{m\times n}[[0,1]^2]$, 
	we define an associated multiplier operator $\text{M}[N]:\R\to L_2^{m}[[0,1]^2]$ and integral operator $\smallint_{[0,1]^2}[N]:L_2^{n}[[0,1]^2]\to\R^{m}$ by
	\begin{align*}
	&(\text{M}[N]u)(x,y):=N(x,y)\mbf{u},	&	&u\in\R^{n}, \\
		&\smallint_{[0,1]^2}[N]\mbf{v}\ :=\int_{0}^{1}\int_{0}^{1}N(\theta,\eta)\mbf{v}(\theta,\eta)d\eta d\theta,	&	&\mbf{v}\in L_2^{n}[[0,1]^2].
	\end{align*}

	\subsection{Algebras of PI Operators on 2D Functions}\label{sec:subsec:2D_PI}	
	
	Partial Integral (PI) operators are bounded, linear operators, parameterized by square integrable functions. We briefly recall the definition and properties of a class of such operators on $L_2[[0,1]^2]$ here, referring to~\cite{jagt2021PIEArxiv} for more details and proofs.
	
	\begin{defn}[2D PI Operators, $\PIset_{2D}$]\label{defn:2D_PI}
		For given parameters $R:=\smallbmat{R_{00}&R_{01}&R_{02}\\R_{10}&R_{11}&R_{12}\\R_{20}&R_{21}&R_{22}}$ with $R_{00}\in L_{2}^{m\times n}[[0,1]^2]$, $R_{0j},R_{i0}\in L_{2}^{m\times n}[[0,1]^3]$ and $R_{ij}\in L_2^{m\times n}[[0,1]^4]$ for $i,j\in\{1,2\}$, define an associated operator $\PI[R]:L_2^{n}[[0,1]^2]\rightarrow L_2^{m}[[0,1]^2]$ such that, for any $\mbf{v}\in L_2^{n}[[0,1]^2]$,
		{
			\begin{align*}
				&(\PI[R]\mbf{v})(x,y):=R_{00}(x,y)\mbf{v}(x,y)\\
				& +\int_{0}^{x}R_{10}(x,y,\theta)\mbf{v}(\theta,y)d\theta +\int_{x}^{1}R_{20}(x,y,\theta)\mbf{v}(\theta,y)d\theta \notag\\[-0.2em] &\quad+\int_{0}^{y}R_{01}(x,y,\eta)\mbf{v}(x,\eta)d\eta +\int_{y}^{1}R_{02}(x,y,\eta)\mbf{v}(x,\eta)d\eta \notag\\[-0.2em]
				&\qquad+\int_{0}^{x}\int_{0}^{y}R_{11}(x,y,\theta,\eta)\mbf{v}(\theta,\eta)d\eta d\theta \notag\\[-0.2em]
				&\qquad\qquad+\int_{x}^{1}\int_{0}^{y}R_{21}(x,y,\theta,\eta)\mbf{v}(\theta,\eta)d\eta d\theta \notag\\[-0.2em]
				&\qquad\qquad\qquad+\int_{0}^{x}\int_{y}^{1}R_{12}(x,y,\theta,\eta)\mbf{v}(\theta,\eta)d\eta d\theta \notag\\[-0.2em]
				&\qquad\qquad\qquad\qquad+\int_{x}^{1}\int_{y}^{1}R_{22}(x,y,\theta,\eta)\mbf{v}(\theta,\eta)d\eta d\theta.	\notag
			\end{align*}
		}
		We refer to an operator $\mcl{R}=\PI[R]$ of this form as a 2D PI operator, writing $\mcl{R}\in \PIset_{2D}^{m\times n}$.
	\end{defn}

	Through slight abuse of notation, we will also use the structure of 2D PI operators to represent maps between $L_{2}[[0,1]^2]$ and the lower-dimensional spaces $L_{2}[0,1]$ and $\R$.
	For example, for $Q:=\smallbmat{Q_{0}&Q_{1}&Q_{2}\\0&0&0\\0&0&0}$ and $R:=\smallbmat{0&0&0\\R_{0}&R_{1}&R_{2}\\R_{0}&R_{1}&R_{2}}$ with $Q_{0},R_{0}\in L_{2}^{m\times n}[[0,1]^2]$ and $Q_{i},R_{i}\in L_2^{m\times n}[[0,1]^3]$ for $i\in\{1,2\}$, the associated PI operators $\PI[Q]:L_2^{n}[0,1]\to L_2^{m}[[0,1]^2]$ and $\PI[R]:L_2^{n}[[0,1]^2]\to L_2^{m}[0,1]$ take the form
	\begin{align*}
		&\bl(\PI[Q]\mbf{u}\br)(x,y)
		=Q_{0}(x,y)\mbf{u}(x)\hspace*{2.1cm}\forall\mbf{u}\in L_2^{n}[0,1], \\
		&\quad +\int_{0}^{x}Q_{1}(x,y,\theta)\mbf{u}(\theta)d\theta  +\int_{x}^{1}Q_{2}(x,y,\theta)\mbf{u}(\theta)d\theta, \\
		&\bl(\PI[R]\mbf{v}\br)(x)=\int_{0}^{1}\bbbl[R_{0}(x,\eta)\mbf{v}(x,\eta)
		\hspace*{1.0cm}\forall\mbf{v}\in L_2^{n}[[0,1]^2], \\[-0.1em]
		&+\int_{0}^{x}R_{1}(x,\theta,\eta)\mbf{v}(\theta,\eta)d\theta +\int_{x}^{1}R_{2}(x,\theta,\eta)\mbf{v}(\theta,\eta)d\theta\bbbr]d\eta.
	\end{align*}
	Similarly, letting $Q:=\smallbmat{Q_{0}&0&0\\Q_{1}&0&0\\Q_{2}&0&0}$ and $R:=\smallbmat{0&R_{0}&R_{0}\\0&R_{1}&R_{1}\\0&R_{2}&R_{2}}$, the associated PI operators $\PI[Q]:L_2^{n}[0,1]\to L_2^{m}[[0,1]^2]$ and $\PI[R]:L_2^{n}[[0,1]^2]\to L_2^{m}[0,1]$ take the form
	\begin{align*}
		&\bl(\PI[Q]\mbf{u}\br)(x,y)
		=Q_{0}(x,y)\mbf{u}(y)\hspace*{2.1cm}\forall\mbf{u}\in L_2^{n}[0,1], \\
		&\quad +\int_{0}^{y}Q_{1}(x,y,\eta)\mbf{u}(\eta)d\eta  +\int_{y}^{1}Q_{2}(x,y,\eta)\mbf{u}(\theta)d\eta, \\
		&\bl(\PI[R]\mbf{v}\br)(y)=\int_{0}^{1}\bbbl[R_{0}(y,\theta)\mbf{v}(\theta,y)
		\hspace*{1.0cm}\forall\mbf{v}\in L_2^{n}[[0,1]^2], \\[-0.1em]
		&+\int_{0}^{y}R_{1}(y,\theta,\eta)\mbf{v}(\theta,\eta)d\eta +\int_{y}^{1}R_{2}(y,\theta,\eta)\mbf{v}(\theta,\eta)d\eta\bbbr]d\theta.
	\end{align*}
	In both cases, we write $\mcl{Q}\in\PIset_{2D\leftarrow 1D}^{m\times n}$ and $\mcl{R}\in\PIset_{1D\leftarrow 2D}^{m\times n}$.
	Finally, for $K\in L_2^{m\times n}[[0,1]^2]$, we define PI operators $\PI\smallbmat{K&0&0\\0&0&0\\0&0&0}=\text{M}[K]:\R^{n}\to L_2^{m}[[0,1]^2]$ and $\PI\smallbmat{0&0&0\\0&K&K\\0&K&K}=\smallint_{[0,1]^2}[K]:L_2^{n}[[0,1]^2]\to\R^{m}$, writing $\tnf{M}[K]\in \PIset_{2D\leftarrow 0}^{m\times n}$ and $\smallint_{[0,1]^2}[K]\in\PIset_{0\leftarrow 2D}^{m\times n}$.

	Throughout this document, we will repeatedly make use of the fact that PI operators form a $*$-algebra, so that the sum $\mcl{Q}+\mcl{R}$, composition $\mcl{Q}\mcl{R}$ and adjoint $\mcl{Q}^{*}$ of any PI operators of suitable dimensions is again a PI operator.
	A proof of this fact, as well as an explicit expression for the parameters defining each of these operation can be found in~\cite{jagt2021PIEArxiv}. Moreover, these operations have also been implemented in the Matlab software suite PIETOOLS (see~\cite{shivakumar2021PIETOOLS}), allowing 2D PI operators $\mcl{Q},\mcl{R}$ to be declared as so-called \texttt{opvar2d} objects \texttt{Q},\texttt{R}, at which point the sum, composition, and adjoint can be readily computed as \texttt{Q+R}, \texttt{Q*R} and \texttt{Q'}, respectively.

\subsection{A PIE Representation of 2D Input-Output PDEs}

A Partial Integral Equation (PIE) is a linear differential equation, parameterized by PI operators, defining the dynamics of a state $\mbf{v}\in L_{2}^{n_{u}}$. For a system with finite-dimensional disturbance $w(t)\in\R^{n_{w}}$ and output $z(t)\in\R^{n_{z}}$, such a PIE takes the form
\begin{align}\label{eq:PIE_standard}
	\mcl{T}\mbf{v}_{t}(t)&=\mcl{A}\mbf{v}(t)+\mcl{B}w(t),	&
	z(t)&=\mcl{C}\mbf{v}(t) +\mcl{D}w(t),
\end{align}
where the parameters ($\mcl{T}$, $\mcl{A}$, etc.) are PI operators. We define solutions to the PIE as follows.
\begin{defn}[Solution to the PIE]\label{defn:PIE_solution}
	For a given input signal $w$, and given initial values $\mbf{v}_{0}\in L_2^{n_{u}}$, we say that $(\mbf{v},z)$ is a solution to the PIE defined by $\{\mcl{T},\mcl{A},\mcl{B},\mcl{C},\mcl{D}\}$ if $\mbf{v}$ is Frech\'et differentiable, $\mbf{v}(0)=\mbf{v}_{0}$, and for all $t\geq 0$, $(\mbf{v}(t),z(t),w(t))$ satisfies~\eqref{eq:PIE_standard}.
\end{defn}
It has previously been shown how we can construct a PIE representation for a broad class of linear 2D PDEs with finite-dimensional inputs and outputs.
In this paper, we focus on 2nd order, coupled, 2D PDEs of the form
\begin{align}\label{eq:PDE_standard}
	\mbf{u}_{t}(t)&=\textstyle{\sum_{i,j=0}^{2}}\text{M}[A_{ij}]\partial_{x}^{i}\partial_{y}^{j}\mbf{u}(t) +\tnf{M}[B]w(t),\hspace*{0.3cm} \mbf{u}(t)\in X, \notag\\
	z(t)&=\textstyle{\sum_{i,j=0}^{2}}\smallint_{[0,1]^2}[C_{ij}]\partial_{x}^{i}\partial_{y}^{j}\mbf{u}(t) +\tnf{M}[D]w(t), 
\end{align}
parameterized by matrix-valued functions
\begin{align*}
	\bmat{A_{ij}&B\\ C_{ij}&D}\in
	\bmat{L_{\infty}^{n_{u}\times n_{u}}[[0,1]^2]&L_2^{n_{u}\times n_{w}}[[0,1]^2]\\
		L_2^{n_{z}\times n_{u}}[[0,1]^2]&\R^{n_{z}\times n_{w}}},
\end{align*}
and where $\mbf{u}(t)\in X\subseteq W_{2}^{n_{u}}[[0,1]^2]$ is the state, $w(t)\in\R^{n_{w}}$ the disturbance, and $z(t)\in \R^{n_{z}}$ the regulated output.
Here, the domain $X$ of the state is assumed to be defined by a set of suitable linear boundary conditions. These boundary conditions may be expressed in terms of all admissible derivatives of the PDE state along the boundary of the domain, collecting these derivatives using the operator $\Lambda_{\tnf{bf}}:=\smallbmat{\Lambda_{1}\\\Lambda_{2}\\\Lambda_{3}}:L_2^{n_{u}}[[0,1]^2]\to \R^{16 n_{u}}\times L_2^{4n_{u}}[0,1]\times L_2^{4n_{u}}[0,1]$, where
\begin{equation}\label{eq:Lambda}
	\Lambda_1\! :=\!
	{\slbmat{\Delta_1 \\
			\Delta_1\partial_x \\
			\Delta_1\partial_y  \\
			\Delta_1\partial_{x}\partial_{y}}}
	,\enspace
	\Lambda_2 \!:=\!\lbmat{\Delta_2 \partial_x^2  \\
		\Delta_2\partial_x^2\partial_y},
	\enspace
	\Lambda_{3} \!:=\!\lbmat{\Delta_3\partial_y^2   \\
		\Delta_3\partial_x\partial_y^2},
\end{equation}
defined by corner, $y$-boundary, and $x$-boundary operators
\begin{equation}\label{eq:Dirac_ops}
	\Delta_1:={\slbmat{
			\Delta_x^{0}\Delta_y^{0}\\
			\Delta_x^{1}\Delta_y^{0}\\
			\Delta_x^{0}\Delta_y^{1}\\
			\Delta_x^{1}\Delta_y^{1}
	}},\quad \Delta_2:= \bmat{
		\Delta_y^{0}    \\
		\Delta_y^{1}
	},\quad \Delta_3:= \bmat{
		\Delta_x^{0}    \\
		\Delta_x^{1}
	}.
\end{equation}
We then parameterize a general class of linear boundary conditions by a matrix $E\in\R^{8n_{u}\times 24n_{u}}$ as
\begin{equation}\label{eq:Xset}
	\mbf{u}(t)\in X:=\{\mbf{u}\in W_{2}^{n_{u}}[[0,1]^2]\mid \text{M}[E]\Lambda_{\tnf{bf}}\mbf{u}=0\},
\end{equation}
which can be used to represent most common boundary conditions (including e.g. Dirichlet, Neumann, Robin). For example, the Dirichlet-Neumann conditions
\begin{equation*}
	\mbf{u}(t,0,y)=\mbf{u}_{x}(t,1,y)=0,\quad \mbf{u}(t,x,0)=\mbf{u}_{y}(t,x,1)=0,
\end{equation*}
can be equivalently represented as
\begin{align*}
	\mbf{u}(t,0,0)=\mbf{u}_{x}(t,1,0)=\mbf{u}_{y}(t,0,1)=\mbf{u}_{xy}(t,1,1)&=0,	\\
	\mbf{u}_{xx}(t,x,0)=\mbf{u}_{xxy}(t,x,1)&=0,	\\
	\mbf{u}_{yy}(t,0,y)=\mbf{u}_{xyy}(t,1,y)&=0,
\end{align*}
which can be readily formulated as in~\eqref{eq:Xset}.
Although, to reduce notation, disturbances in the boundary conditions will not be considered here, such disturbances can also be included, using the methodology presented in~\cite{jagt2022PIE_2DHinfty}.

Given a PDE state $\mbf{u}(t)\in X$, we define an associated \textit{fundamental state} as its highest-order spatial derivative, $\mbf{v}(t):=\partial_{x}^{2}\partial_{y}^{2}\mbf{u}(t)\in L_{2}$, which is free of the boundary conditions and continuity constraints imposed upon $\mbf{u}(t)$. Then, if the boundary conditions are sufficiently well-posed (excluding e.g. periodic boundary conditions), it has been shown that there exists a PI operator $\mcl{T}\in\PIset_{2D}^{n_{u}\times n_{u}}$ such that $\mbf{u}(t)=\mcl{T}\mbf{v}(t)$. Moreover, in this case, we can define PI operators $\{\mcl{A},\mcl{B},\mcl{C},\mcl{D}\}$ such that $\mbf{u}(t)$ satisfies the PDE~\eqref{eq:PDE_standard} if and only if $\mbf{v}(t)$ satisfies the PIE~\eqref{eq:PIE_standard}.
In particular, we have the following result, corresponding to Lemma~14 in~\cite{jagt2022PIE_2DHinfty}.


\begin{lem}\label{lem:PDE_to_PIE_io}
	Let parameters $\{A_{ij},B,C_{ij},D,E\}$ define a PDE as in~\eqref{eq:PDE_standard} with boundary conditions as in~\eqref{eq:Xset}, and satisfy the conditions in Lemma~14 in~\cite{jagt2022PIE_2DHinfty}. Let PI operators $\{\mcl{T},\mcl{A},\mcl{B},\mcl{C},\mcl{D}\}$ be as defined in that lemma. Then, for any $\mbf{v}\in L_2^{n_{u}}[[0,1]^2]$ and $\mbf{u}\in X$,
	\begin{equation*}
		\mbf{u} = \mcl{T}\:\partial_{x}^2\partial_{y}^2\mbf{u}
		\quad\text{and}\quad
		\mbf{v} = \partial_{x}^2\partial_{y}^2\:\mcl{T}\mbf{v}.
	\end{equation*}
	Moreover, for any input $w(t)\in \R^{n_{w}}$ and time $t\geq 0$, $(\mbf{u}(t),z(t))\in X\times\R^{n_{z}}$ satisfies the PDE~\eqref{eq:PDE_standard} if and only if $(\partial_{x}^2\partial_{y}^2\mbf{u}(t),z(t))\in L_2^{n_{u}}\times\R^{n_{z}}$ satisfies the PIE~\eqref{eq:PIE_standard}. Conversely, $(\mbf{v}(t),z(t))\in L_2^{n_{u}}\times\R^{n_{z}}$ satisfies the PIE~\eqref{eq:PIE_standard} if and only if $(\mcl{T}\mbf{v}(t),z(t))\in X\times\R^{n_{z}}$ satisfies the PDE~\eqref{eq:PDE_standard}.	
\end{lem}

Note that, for a given PDE, the associated PIE representation as defined in Lem.~\ref{lem:PDE_to_PIE_io} is free of any boundary conditions. Instead, these boundary conditions are incorporated into the definitions of the PI operators $\{\mcl{T},...,\mcl{D}\}$ defining the PIE representation, so that in particular $\mcl{T}\mbf{v}(t)\in X$ for any $\mbf{v}(t)\in L_{2}^{n_{u}}$.
We refer to earlier works (e.g.~\cite{jagt2022PIE_2DHinfty}) for more information on computing these PI operators for a given PDE.
This conversion procedure has also been fully incorporated in the PIETOOLS software suite, allowing a 2D PDE as in~\eqref{eq:PDE_standard} (as well as more general linear 1D and 2D PDEs) to be readily declared as a structure \texttt{PDE}, at which point the associated PIE representation can be computed by simply calling \texttt{PIE=convert(PDE)}.

\subsection{A Linear PI Inequality (LPI) for $L_2$-gain Analysis}

Linear PI Inequalities (LPIs) are convex optimization programs involving linear operator inequalities on PI operator variables. Exploiting the algebraic properties of PI operators, and the lack of auxiliary constraints (boundary conditions) in the PIE representation, several problems of analysis and control of PIEs can be posed as LPIs. For example, the following LPI for $L_{2}$-gain analysis of PIEs was derived in Lem.~8 in~\cite{jagt2022PIE_2DHinfty}.
%
\begin{lem}\label{lem:KYP}
	Let $\gamma>0$, and suppose there exists a PI operator $\mcl{P}\in\PIset_{2D}^{n_u\times n_u}$ such that $\mcl{P}=\mcl{P}^*\succ 0$ and
	\begin{align*}
		&\bmat{-\gamma I & \mcl{D} &~ \mcl{C}\\
			\mcl{D}^* & -\gamma I &~ \mcl{B}^*\mcl{P}\mcl{T}	\\
			\mcl{C}^* & \mcl{T}^*\mcl{P}\mcl{B} &~ \mcl{A}^*\mcl{P}\mcl{T} + \mcl{T}^*\mcl{P}\mcl{A}}	
		\preceq 0.
	\end{align*}
	Then, for any $w\in L_2^{n_w}[0,\infty)$, if $(w,z)$ satisfies the PIE~\eqref{eq:PIE_standard} with $\mbf{v}(0)=\mbf{0}$, then $z\in L_2^{n_z}[0,\infty)$ and $\|z\|_{L_2}\leq\gamma \|w\|_{L_2}$.
\end{lem}

We use this result to derive a similar LPI for $H_{\infty}$-optimal estimator synthesis for 2D PDEs in the next section.

\section{An $H_{\infty}$-Optimal Estimator for 2D PDEs}\label{sec:main_result}

In this section, we provide the main technical result of this paper, proposing 
an LPI for $H_{\infty}$-optimal estimator synthesis for a class of 2D PDEs as in~\eqref{eq:PDE_standard}. 
We suppose that we have three observed output signals defined as
\begin{align}\label{eq:PDE_outputs}
	\mbf{q}(t)&:=\slbmat{q_{1}(t)\\\mbf{q}_{2}(t)\\\mbf{q}_{3}(t)}=\slbmat{\text{M}[C_{1}]\Lambda_{1}\\\text{M}[C_{2}]\Lambda_{2}\\\text{M}[C_{3}]\Lambda_{3}}\mbf{u}(t) +\slbmat{\text{M}[D_{1}]\\\text{M}[D_{2}]\\\text{M}[D_{3}]}w(t),
\end{align}
with $q_{1}(t)\in\R^{n_{q_1}}$, $\mbf{q}_{2}(t)\in L_2^{n_{q_{2}}}[0,1]$ and $\mbf{q}_{3}(t)\in L_2^{n_{q_{3}}}[0,1]$, 
defined by parameters
\begin{equation*}
	\bmat{C_{1}&C_{2}&C_{3}\\ D_{1}&D_{2}&D_{3}}\!\in\!\bmat{\R^{n_{q_{1}}\times 16n_{u}} & L_2^{n_{q_{2}}\times 4n_{u}}[0,1] & L_2^{n_{q_{3}}\times 4n_{u}}[0,1]\\ \R^{n_{q_{1}}\times n_{w}} & L_2^{n_{q_{2}}\times n_{w}}[0,1] & L_2^{n_{q_{3}}\times n_{w}}[0,1]},
\end{equation*}
and where the trace operators $\Lambda_{i}$ for $i\in\{1,2,3\}$
are as defined in~\eqref{eq:Lambda}, evaluating admissible derivatives of the state along the boundary of the domain. 
We define a solution to the resulting system as follows.
\begin{defn}[Solution to the PDE]\label{defn:PDE_solution}
	For a given input signal $w$ and initial state $\mbf{u}_{0}\in X$, we say that $(\mbf{u},z,\mbf{q})$ is a solution to the PDE defined by $\{A_{ij},B,C_{ij},D,C_{k},D_{k},E\}$ if $\mbf{u}$ is Frech\'et differentiable, $\mbf{u}(0)=\mbf{u}_{0}$, and for all $t\geq 0$, $(\mbf{u}(t),z(t),\mbf{q}(t))$ satisfies~\eqref{eq:PDE_standard} and~\eqref{eq:PDE_outputs}.
\end{defn}

Now, to construct an estimator for the PDE with the proposed output, first note that (by Lem.~\ref{lem:PDE_to_PIE_io}) the dynamics of the PDE can be equivalently represented in terms of the fundamental state $\mbf{v}(t)=\partial_{x}^{2}\partial_{y}^{2}\mbf{u}(t)$, as the PIE~\eqref{eq:PIE_standard}. Invoking the identity $\mbf{u}=\mcl{T}\mbf{v}$, the output signals can also be represented in terms of this fundamental state as
\begin{equation}\label{eq:PIE_outputs}
	\mbf{q}(t)=\mcl{C}_{q}\mbf{v}(t) +\mcl{D}_{q}w(t),
\end{equation}
where we define the operators 
\begin{equation}\label{eq:observer_ops}
	\mcl{C}_{q}:={\slbmat{\mcl{C}_{1}\\\mcl{C}_{2}\\\mcl{C}_{3}}}
	={\slbmat{\text{M}[C_{1}]\circ\Lambda_{1}\circ\mcl{T}\\\text{M}[C_{2}]\circ\Lambda_{2}\circ\mcl{T}\\\text{M}[C_{2}]\circ\Lambda_{2}\circ\mcl{T}}},\quad
	\mcl{D}_{q}:={\slbmat{\text{M}[D_{1}]\\\text{M}[D_{2}]\\\text{M}[D_{3}]}}.
\end{equation}
Then, a Luenberger-type estimator for the PIE~\eqref{eq:PIE_standard} can be parameterized by a PI operator $\mcl{L}$ as
\begin{align}\label{eq:PIE_estimator}
	\mcl{T}\hat{\mbf{v}}_{t}(t)&=\!\mcl{A}\hat{\mbf{v}}(t) +\mcl{L}(\mcl{C}_{q}\hat{\mbf{v}}(t)-\mbf{q}(t)),	&
	\hat{z}(t)&=\mcl{C}\hat{\mbf{v}}(t),
\end{align}
returning an estimate of the PDE state $\mbf{u}(t)=\mcl{T}\mbf{v}(t)$ as $\hat{\mbf{u}}(t)=\mcl{T}\hat{\mbf{v}}(t)$. 
The goal, then, is to choose the gain $\mcl{L}$ such as to minimize the $H_{\infty}$-norm of the resulting error dynamics, i.e. to solve the optimization program
\begin{align*}
	\min_{\mcl{L},\gamma}~\gamma\hspace*{0.4cm}
	\tnf{s.t.}\quad \|\hat{z}-z\|_{L_{2}}\leq \gamma\|w\|_{L_{2}},~\forall w\in L_{2}[0,\infty)\setminus\{0\}.
\end{align*}
The following result shows that a solution to this program can be computed by solving an LPI.

\begin{thm}\label{thm:PDE_estimator}
	For given $\tnf{G}_{\tnf{pde}}:=\{A_{ij},B,C_{ij},D,C_{k},D_{k},E\}$, define associated PI operators $\{\mcl{T},\mcl{A},\mcl{B},\mcl{C},\mcl{D}\}$ as in Lem.~\ref{lem:PDE_to_PIE_io}, and let further $\{\mcl{C}_{q},\mcl{D}_{q}\}$ be as in~\eqref{eq:observer_ops}.	
	Suppose that there exists a constant $\gamma>0$ and PI operators $\mcl{P}\in\PIset_{2D}^{n_u\times n_u}$ and $\mcl{W}:\PIset_{2D\leftarrow 0}^{n_{u}\times n_{q_1}}\times \PIset_{2D\leftarrow 1D}^{n_{u}\times n_{q_{2}}}\times \PIset_{2D\leftarrow 1D}^{n_{u}\times n_{q_3}}$ such that \\[-0.55em]
	\begin{equation*}
		\mcl{P}=\mcl{P}^*\succ0,	\quad
		\bmat{-\gamma I &\! -\mcl{D} &\! \mcl{C}\\
			(\cdot)^* &\! -\gamma I &\! -[\mcl{B}^*\mcl{P}+\mcl{D}_{q}^*\mcl{W}^*]\mcl{T}	\\
			(\cdot)^* &\! (\cdot)^* &\! (\cdot)^* + \mcl{T}^*[\mcl{P}\mcl{A}+\mcl{W}\mcl{C}_{q}]}	
		\preceq 0,
	\end{equation*}
	and let $\mcl{L}=\mcl{P}^{-1}\mcl{W}$.
	Then, if $(\mbf{u},z,\mbf{q})$ is a solution to the PDE defined by $\tnf{G}_{\tnf{pde}}$ for some disturbance $w\in L_2^{n_w}[0,\infty)$ and initial state $\mbf{u}_{0}\in X$, and $(\hat{\mbf{v}},\hat{z})$ is a solution to the PIE~\eqref{eq:PIE_estimator} with input $\mbf{q}$ and initial state $\hat{\mbf{v}}_{0}=\partial_{x}^2\partial_{y}^{2}\mbf{u}_{0}$, then $\hat{z}-z\in L_2^{n_z}[0,\infty)$ and $\|\hat{z}-z\|_{L_2}\leq\gamma \|w\|_{L_2}$.
\end{thm}

To prove this result, in Subsec.~\ref{sec:main_result:PIE}, it is first proven that the operator $\mcl{C}_{q}$ in~\eqref{eq:PIE_outputs} is indeed a PI operator, thus yielding a PIE representation of the considered PDE.
In Subsec.~\ref{sec:main_result:LPI}, it is then shown how an estimator for this PIE can be synthesized by solving the proposed LPI.

\subsection{Representation of Infinite-Dimensional PDE Outputs}\label{sec:main_result:PIE}

In order to synthesize an estimator for the PDE~\eqref{eq:PDE_standard} with observed outputs as in~\eqref{eq:PDE_outputs}, we will first represent the system with these outputs as an equivalent PIE. Here, Lem.~\ref{lem:PDE_to_PIE_io} already shows that we can define PI operators $\{\mcl{T},\mcl{A},\mcl{B},\mcl{C},\mcl{D}\}$ to express the PDE dynamics as a PIE~\eqref{eq:PIE_standard}, modeling the evolution of the fundamental state $\mbf{v}(t)=\partial_{x}^2\partial_{y}^{2}\mbf{u}(t)\in L_2^{n_{u}}$ associated to the PDE state $\mbf{u}(t)\in X$. Moreover, substituting the relation $\mbf{u}(t)=\mcl{T}\mbf{v}(t)$ into the expression for the observed outputs in~\eqref{eq:PDE_outputs}, it is clear that these outputs can be equivalently expressed as
\begin{align*}
	\mbf{q}(t)&:={\slbmat{\text{M}[C_{1}]\Lambda_{1}\\\text{M}[C_{2}]\Lambda_{2}\\\text{M}[C_{3}]\Lambda_{3}}}\mcl{T}\mbf{v}(t) +{\slbmat{\text{M}[D_{1}]\\\text{M}[D_{2}]\\\text{M}[D_{3}]}}w(t).
\end{align*}
It remains only to prove, then, that the operators $\Lambda_{i}\circ\mcl{T}$ for $\Lambda_{i}$ as in~\eqref{eq:Lambda} are indeed PI operators. 
To begin, consider the first element of the operator $\Lambda_{2}$, defined as $\Delta_{y}^{0}\partial_{x}^{2}$. 
In order to prove that the composition $\Lambda_{2}\circ\mcl{T}$ is a PI operator, we need to show that for all $\mbf{v}\in L_2^{n_{u}}$, we can express $\partial_{x}^2(\mcl{T}\mbf{v})(x,y)=(\mcl{R}\mbf{v})(x,y)$ for some PI operator $\mcl{R}$, and subsequently that $(\Delta_{y}^{0}\mcl{R}\mbf{v})(x)=(\mcl{Q}\mbf{v})(x)$ for some PI operator $\mcl{Q}$. For this first part, defining a PI operator $\mcl{R}=\partial_{x}^2\:\mcl{T}$, we remark that composition rules for partial differential operators and PI operators have already been derived in~\cite{jagt2021PIEArxiv}, and in particular, we have the following result.

\begin{lem}\label{lem:PI_derivative}
	For given $E\in\R^{8n_{u}\times 24n_{u}}$ defining a set $X\subseteq W_{2}^{n_{u}}[[0,1]^2]$ as in~\eqref{eq:Xset}, let $T:=\smallbmat{0&0&0\\0&T_{11}&T_{12}\\0&T_{21}&T_{22}}$ for $T_{ij}\in L_2^{n_{u}\times n_{u}}[[0,1]^2]$ be the associated parameters as defined in Thm.~13 in~\cite{jagt2021PIEArxiv}, so that $\mcl{T}=\PI[T]\in\PIset_{2D}^{n_{u}\times n_{u}}$ is as in Lem~\ref{lem:PDE_to_PIE_io}. For $k,\ell\in\{0,1,2\}$, define parameters $R^{k\ell}:=\smallbmat{R_{00}^{k\ell}&R_{01}^{k\ell}&R_{02}^{k\ell}\\R_{10}^{k\ell}&R_{11}^{k\ell}&R_{12}^{k\ell}\\R_{20}^{k\ell}&R_{21}^{k\ell}&R_{22}^{k\ell}}$ by
	\begin{align*}
		&R_{ij}^{k\ell}(x,y,\theta,\eta):=\bbbl\{\slmat{\partial_{x}^{k}\partial_{y}^{\ell}T_{ij}(x,y,\theta,\eta),\quad	k<2,~\ell<2, \\
			0,\hspace*{2.5cm}	\text{else},}\\
		&R_{0j}^{k\ell}(x,y,\eta):=\bbbl\{\slmat{R_{1j}^{1\ell}(x,y,x,\eta)-R_{2j}^{1\ell}(x,y,x,\eta),\quad	k=2,~\ell<2,	\\
			0,\hspace*{4.1cm}	\text{else},}\\
		&R_{i0}^{k\ell}(x,y,\theta):=\bbbl\{\slmat{R_{i1}^{k1}(x,y,\theta,y)-R_{i2}^{k1}(x,y,\theta,y),\quad	k<2,~\ell=2,	\\
			0,\hspace*{4.1cm}	\text{else},}\\
		&R_{00}^{k\ell}(x,y):=\bbbl\{\slmat{I_{n_{u}},\quad	k=2,~\ell=2,	\\
			0,\qquad	\text{else},}
	\end{align*}
	for $i,j\in\{1,2\}$. Then, for every $\mbf{v}\in L_2^{n_{u}}$,
	\begin{equation*}
		\partial_{x}^{k}\partial_{y}^{\ell}\bl(\PI[T]\mbf{v}\br)(x,y)
		=\bl(\PI[R^{k\ell}]\mbf{v}\br)(x,y),\qquad \forall k,\ell\in\{0,1,2\}.
	\end{equation*}
\end{lem}
\begin{pf}
	The result follows by the Leibniz integral rule and the definition of the paameters $T_{ij}$. An explicit derivation is given in the proof of Thm.~13 in~\cite{jagt2021PIEArxiv}.
\end{pf}

Lem.~\ref{lem:PI_derivative} shows that, for each $0\leq k,\ell\leq 2$, we can explicitly define 2D PI operators $\mcl{R}_{k\ell}\in\PIset_{2D}^{n_{u}\times n_{u}}$ such that $\partial_{x}^{k}\partial_{y}^{\ell}\circ\mcl{T}=\mcl{R}_{k\ell}$. 
By definition of the operators $\Lambda_{i}$ in~\eqref{eq:Lambda}, then, it follows that
\begin{equation*}
	\Lambda_1\mcl{T} =
	{\slbmat{\Delta_1\mcl{T} \\
		\Delta_1\mcl{R}_{10} \\
		\Delta_1\mcl{R}_{01}  \\
		\Delta_1\mcl{R}_{11}}}
	,\enspace
	\Lambda_2\mcl{T} =\lbmat{\Delta_2 \mcl{R}_{20}  \\
		\Delta_2 \mcl{R}_{21}},
	\enspace
	\Lambda_{3}\mcl{T} =\lbmat{\Delta_3\mcl{R}_{02}   \\
		\Delta_3\mcl{R}_{12}},
\end{equation*}
where the Dirac operators $\Delta_{1}$, $\Delta_{2}$, and $\Delta_{3}$ are as in~\eqref{eq:Dirac_ops}, evaluating the state at the corners, $y$-boundaries, and $x$-boundaries of the domain, respectively. It remains to show that the composition of these Dirac operators with the different PI operators $\mcl{R}_{k\ell}$ can also be expressed as PI operators.
For this, we remark that e.g. evaluating the partial integral $\int_{0}^{x}R(x,y,\theta)\mbf{v}(\theta,y)d\theta$ at $x=1$, we can express the result as a full integral operator on $\mbf{v}$ as $\int_{0}^{1}R(1,y,\theta)\mbf{v}(\theta,y)d\theta$. The following proposition generalizes this result to compositions of more general Dirac operators with 2D PI operators.

\begin{prop}\label{prop:PI_Dirac}
	Let $R_{i0},R_{0j}\in L_{2}^{m\times n}[[0,1]^2\times[0,1]]$ and $R_{ij}\in L_{2}^{m\times n}[[0,1]^2\times[0,1]^2]$ for $i,j\in\{1,2\}$, and define $F_{j}^{k},G_{i}^{k}\in L_2^{m\times n}[[0,1]\times[0,1]^2]$ and $F_{0}^{k},G_{0}^{\ell},H^{k\ell}\in L_2^{m\times n}[[0,1]^2]$ for $k\in\{0,1\}$ by
	\begin{align*}
		&F_{0}^{k}(y,\theta):=\bbbl\{\slmat{R_{20}(0,y,\theta),\quad k=0,	\\[0.2em]
			R_{10}(1,y,\theta),\quad k=1,} 	\\
		&F_{j}^{k}(y,\theta,\eta):=\bbbl\{\slmat{R_{2j}(0,y,\theta,\eta),\quad k=0,	\\[0.2em]
											R_{1j}(1,y,\theta,\eta),\quad k=1,}	\\
		&G_{0}^{\ell}(x,\eta):=\bbbl\{\slmat{R_{02}(x,0,\eta),\quad \ell=0,	\\[0.2em]
			R_{01}(x,1,\eta),\quad \ell=1,}		\\
		&G_{i}^{\ell}(x,\theta,\eta):=\bbbl\{\slmat{R_{i2}(x,0,\theta,\eta),\quad \ell=0,\\[0.2em]
												R_{i1}(x,1,\theta,\eta),\quad \ell=1,}	\\
		&H^{k\ell}(\theta,\eta):=\bbbl\{\slmat{F_{2}^{k}(0,\theta,\eta),\quad \ell=0,\\[0.2em]
			F_{1}^{k}(1,\theta,\eta),\quad \ell=1,}
	\end{align*}
	for $x,y,\theta,\eta\in[0,1]$. Then
	\begin{align*}
		&\Delta_{x}^{k}\PI\smallbmat{0&0&0\\R_{10}&R_{11}&R_{12}\\R_{20}&R_{21}&R_{22}}
		=\PI\smallbmat{0&0&0\\F_{0}^{k}&F_{1}^{k}&F_{2}^{k}\\F_{0}^{k}&F_{1}^{k}&F_{2}^{k}}
		&	&\hspace*{-0.4cm}\in\PIset_{1D\from 2D}^{m\times n},\\
		&\Delta_{y}^{\ell}\PI\smallbmat{0&R_{01}&R_{02}\\0&R_{11}&R_{12}\\0&R_{21}&R_{22}}
		=\PI\smallbmat{0&G_{0}^{\ell}&G_{0}^{\ell}\\0&G_{1}^{\ell}&G_{1}^{\ell}\\0&G_{2}^{\ell}&G_{2}^{\ell}}
		&	&\hspace*{-0.4cm}\in\PIset_{1D\from 2D}^{m\times n},	\\
		&\Delta_{x}^{k}\Delta_{y}^{\ell}\PI\smallbmat{0&0&0\\0&R_{11}&R_{12}\\0&R_{21}&R_{22}}
		=\PI\smallbmat{0&0&0\\0&H^{k\ell}&H^{k\ell}\\0&H^{k\ell}&H^{k\ell}}	&	&\hspace*{-0.4cm}\in \PIset_{0\from 2D}^{m\times n}.
	\end{align*}
\end{prop}
\begin{pf}
	We prove the result only for $k=\ell=0$, as the proof for each other combination of $k,\ell\in\{0,1\}$ is similar.
	Fix arbitrary $\mbf{v}\in L_2^{n}[[0,1]^2]$.
	Then, evaluating $\bl(\PI\smallbmat{0&0&0\\R_{10}&R_{11}&R_{12}\\R_{20}&R_{21}&R_{22}}\mbf{v}\br)(x,y)$ at $x=0$, the integral terms $\int_{0}^{x}[R_{1j}]$ vanish, and we find
	\begin{align*}
		&\left(\Delta_{x}^{0}\PI\smallbmat{0&0&0\\R_{10}&R_{11}&R_{12}\\R_{20}&R_{21}&R_{22}}\mbf{v}\right)(y)	
		=\int_{0}^{1}\bbbl[R_{20}(0,y,\theta)\mbf{v}(\theta,y) \\ &+\!\!\int_{0}^{y}\!\!R_{21}(0,y,\theta,\eta)\mbf{v}(\theta,\eta)d\eta \!+\!\!\int_{y}^{1}\!\!R_{22}(0,y,\theta,\eta)\mbf{v}(\theta,\eta)d\eta\bbbr]d\theta	\\[-0.1em]
		&\hspace*{0.0cm}
		\bbl(\PI\smallbmat{0&0&0\\R_{20}(0,.)&R_{21}(0,.)&R_{22}(0,.)\\R_{20}(0,.)&R_{21}(0,.)&R_{22}(0,.)}\mbf{v}\bbr)(y)
		=\bbl(\PI\smallbmat{0&0&0\\F_{0}^{0}&F_{1}^{0}&F_{2}^{0}\\F_{0}^{0}&F_{1}^{0}&F_{2}^{0}}\mbf{v}\bbr)(y).
	\end{align*}
	Similarly, evaluating $\bl(\PI\smallbmat{0&R_{01}&R_{02}\\0&R_{11}&R_{12}\\0&R_{21}&R_{22}}\mbf{v}\br)(x,y)$ at $y=0$, the integral terms $\int_{0}^{y}[R_{i1}]$ vanish, and we obtain
	\begin{align*}
		&\left(\Delta_{y}^{0}\PI\smallbmat{0&R_{01}&R_{02}\\0&R_{11}&R_{12}\\0&R_{21}&R_{22}}\mbf{v}\right)(y)	
		=\int_{0}^{1}\bbbl[R_{02}(x,0,\eta)\mbf{v}(x,\eta) \\ &+\!\!\int_{0}^{x}\!\!R_{12}(x,0,\theta,\eta)\mbf{v}(\theta,\eta)d\theta \!+\!\!\int_{x}^{1}\!\!R_{22}(x,0,\theta,\eta)\mbf{v}(\theta,\eta)d\theta\bbbr]d\eta	\\[-0.1em]
		&\hspace*{5.0cm}=\bbl(\PI\smallbmat{0&G_{0}^{0}&G_{0}^{0}\\0&G_{1}^{0}&G_{1}^{0}\\0&G_{2}^{0}&G_{2}^{0}}\mbf{v}\bbr)(x).
	\end{align*}
	Finally, evaluating $\bl(\PI\smallbmat{0&0&0\\0&R_{11}&R_{12}\\0&R_{21}&R_{22}}\mbf{v}\br)(x,y)$ at $x=y=0$, the integrals $\int_{0}^{x}\!\int_{0}^{y}[R_{11}]$, $\int_{0}^{x}\!\int_{y}^{1}[R_{12}]$, and $\int_{x}^{1}\!\int_{0}^{y}[R_{12}]$ vanish, and we find
	\begin{align*}
		&\left(\Delta_{x}^{0}\Delta_{y}^{0}\PI\smallbmat{0&0&0\\0&R_{11}&R_{12}\\0&R_{21}&R_{22}}\mbf{v}\right)
		=\int_{0}^{1}\!\int_{0}^{1}\!\!R_{22}(0,0,\theta,\eta)\mbf{v}(\theta,\eta)d\eta d\theta \\[-0.1em]
		&\hspace*{5.0cm}
		=\bbl(\PI\smallbmat{0&0&0\\0&H^{00}&H^{00}\\0&H^{00}&H^{00}}\mbf{v}\bbr).
	\end{align*}
\end{pf}

Prop.~\ref{prop:PI_Dirac} proves that, given an operator $\PI[R]\in\PIset_{2D}$ defined by parameters 
with a suitable structure, we can define a PI operator $\PI[G]\in\PIset_{2D}$ such that e.g.
\begin{equation*}
	\bl(\PI[R]\mbf{v}\br)(x,0)=\bl(\Delta_{y}^{0}\PI[R]\mbf{v}\br)(x)=\bl(\PI[G]\mbf{v}\br)(x)
\end{equation*}
Returning now to e.g. the composition $\Delta_{y}^{0}\partial_{x}^2\:\mcl{T}=\Delta_{y}^{0}\mcl{R}_{20}$ where $\mcl{R}_{20}:=\partial_{x}^{2}\:\mcl{T}$ is as in Lem.~\ref{lem:PI_derivative}, the operator $\mcl{R}_{20}$ indeed has a suitable structure that allows us to define a 2D PI operator $\mcl{G}_{20}:=\Delta_{y}^{0}\mcl{R}_{20}$ as in Prop.~\ref{prop:PI_Dirac}.
More generally, we find that each of the compositions $\Lambda_{i}\circ\mcl{T}$ -- and therefore the operator $\mcl{C}_{q}$ in~\eqref{eq:observer_ops} -- can be expressed as a 2D PI operator. Therefore, we can express the infinite-dimensional output $\mbf{q}(t)$ in terms of the fundamental state $\mbf{v}(t)=\partial_{x}^{2}\partial_{y}^{2}\mbf{u}(t)$ using PI operators as in~\eqref{eq:PIE_outputs}, and we obtain a PIE representation of the PDE with sensed outputs.

\begin{lem}\label{lem:PDE2PIE_observed}
	For given $\tnf{G}_{\tnf{pde}}:=\{A_{ij},B,C_{ij},D,C_{k},D_{k},E\}$, define associated PI operators $\tnf{G}_{\tnf{pie}}\!:=\{\mcl{T}\!,\mcl{A},\mcl{B},\smallbmat{\mcl{C}\\\mcl{C}_{q}}\!,\smallbmat{\mcl{D}\\\mcl{D}_{q}}\}$ as in Lem.~\ref{lem:PDE_to_PIE_io} and in~\eqref{eq:observer_ops}. Then, for any input $w$, $(\mbf{u},z,\mbf{q})$ is a solution to the PDE defined by $\tnf{G}_{\tnf{pde}}$ with initial state $\mbf{u}_{0}\in X$ if and only if $(\mbf{v},\smallbmat{z\\\mbf{q}})$ with $\mbf{v}=\partial_{x}^{2}\partial_{y}^{2}\mbf{u}$ is a solution to the PIE defined by $\tnf{G}_{\tnf{pie}}$ with initial state $\mbf{v}_{0}=\partial_{x}^2\partial_{y}^2\mbf{u}_{0}$.
\end{lem}
\begin{pf}
	Let operators $\tnf{G}_{\tnf{pie}}$ be as defined, and fix an arbitrary input $w$. Then, by Lem~\ref{lem:PDE_to_PIE_io}, $(\mbf{u},z,\mbf{q})$ is a solution to the PDE defined by $\tnf{G}_{\tnf{pde}}$ with initial state $\mbf{u}_{0}\in X$ if and only if $(\mbf{v},z)$ is a solution to the PIE defined by $\{\mcl{T},\mcl{A},\mcl{B},\mcl{C},\mcl{D}\}$ with initial state $\mbf{v}_{0}=\partial_{x}^2\partial_{y}^2\mbf{u}_{0}$, and $\mbf{q}(t)$ satisfies~\eqref{eq:PDE_outputs} with $\mbf{u}(t)=\mcl{T}\mbf{v}(t)$. Here, by definition of the operators $\{\mcl{C}_{q},\mcl{D}_{q}\}$, $\mbf{q}(t)$ satisfies~\eqref{eq:PDE_outputs} with $\mbf{u}(t)=\mcl{T}\mbf{v}(t)$ if and only if $\mbf{q}(t)$ satisfies~\eqref{eq:PIE_outputs}, and hence $(\mbf{v},\smallbmat{z\\\mbf{q}})$ is a solution to the PIE defined by $\tnf{G}_{\tnf{pie}}$.
\end{pf}

\subsection{An LPI for Optimal Estimation of PIEs}\label{sec:main_result:LPI}

Having derived a PIE representation of the 2D PDE~\eqref{eq:PDE_standard}, consider now a Luenberger-type estimator for this PIE as in~\eqref{eq:PIE_estimator}, parameterized by a PI operator $\mcl{L}\in \PIset_{2D\leftarrow 0}^{n_{u}\times n_{q_1}}\times \PIset_{2D\leftarrow 1D}^{n_{u}\times n_{q_{2}}}\times \PIset_{2D\leftarrow 1D}^{n_{u}\times n_{q_3}}$.
Then, for any solution $(\mbf{v},z)$ and $(\hat{\mbf{v}},\hat{z})$ to the PIE~\eqref{eq:PIE_standard} and the PIE~\eqref{eq:PIE_estimator}, respectively,
the errors $\mbf{e}(t):=\hat{\mbf{v}}(t)-\mbf{v}(t)$ and $\tilde{z}(t):=\hat{z}(t)-z(t)$ will satisfy
\begin{align}\label{eq:PIE_error}
	\mcl{T}\dot{\mbf{e}}(t)&=\tilde{\mcl{A}}\mbf{e}(t) +\tilde{\mcl{B}}w(t),	&	
	\tilde{z}(t)&=\tilde{\mcl{C}}\mbf{e}(t)+\tilde{\mcl{D}}w(t),
\end{align}
where we define the PI operators
\begin{align*}
	\tilde{\mcl{A}}&:=\mcl{A}+\mcl{L}\mcl{C}_{q}, &
	\tilde{\mcl{B}}&:=-\bl(\mcl{B}+\mcl{L}\mcl{D}_{q}\br),	&
	\tilde{\mcl{C}}&:=\mcl{C},	&
	\tilde{\mcl{D}}&:=-\mcl{D}.
\end{align*}
The $H_{\infty}$-optimal estimator synthesis problem, then, is to establish a value of the operator $\mcl{L}$ that minimizes the $L_2$-gain $\sup_{w\in L_2\setminus\{0\}}\frac{\|\tilde{z}\|_{L_2}}{\|w\|_{L_2}}$ from disturbances $w$ to the error $\tilde{z}$. However, the problem of computing an upper bound $\gamma$ on the $L_{2}$-gain of PIEs with finite-dimensional inputs and outputs has already been tackled in~\cite{jagt2022PIE_2DHinfty}. Using that result, we can readily pose the problem of verifying an upper bound $\gamma$ on the $H_{\infty}$-norm of the error dynamics as follows.
\begin{cor}\label{cor:KYP_estimator}
	Let $\gamma>0$, and suppose there exists PI operators $\mcl{P}\in\PIset_{2D}^{n_{u}\times n_{u}}$ and $\mcl{W}\in\PIset_{2D\leftarrow 0}^{n_{u}\times n_{q_1}}\times \PIset_{2D\leftarrow 1D}^{n_{u}\times n_{q_{2}}+n_{q_{3}}}$ such that
	\begin{align}\label{eq:KYP_estimator_LPI}
		\mcl{P}=\mcl{P}^*\succ0,\enspace 
		&Q:=\bmat{-\gamma I &\! -\mcl{D} &\! \mcl{C}\\
			(\cdot)^* &\! -\gamma I &\! -[\mcl{B}^*\mcl{P}+\mcl{D}_{q}^*\mcl{W}^*]\mcl{T}	\\
			(\cdot)^* &\! (\cdot)^* &\! (\cdot)^* + \mcl{T}^*[\mcl{P}\mcl{A}+\mcl{W}\mcl{C}_{q}]}	
		\preceq 0. \notag\\[-1.0em]
	\end{align}
	Then, for any $w\in L_2^{n_w}[0,\infty)$, if $(w,\tilde{z})$ satisfies the PIE~\eqref{eq:PIE_error} with $\mcl{L}=\mcl{P}^{-1}\mcl{W}$ and $\mbf{e}(0)=\mbf{0}$, then $\tilde{z}\in L_2^{n_z}[0,\infty)$ and $\|\tilde{z}\|_{L_2}\leq\gamma \|w\|_{L_2}$.
\end{cor}
\begin{pf}
	Let the conditions of the corollary be satisfied for some $(\gamma,\mcl{P},\mcl{W})$, and let $\mcl{L}:=\mcl{P}^{-1}\mcl{W}$. Then, by definition of the operators $\{\tilde{\mcl{A}},\tilde{\mcl{B}},\tilde{\mcl{C}},\tilde{\mcl{D}}\}$ defining the PIE~\eqref{eq:PIE_error}, we have
	\begin{equation*}
		[\mcl{B}^*\mcl{P}+\mcl{D}_{q}^*\mcl{W}^*]=-\tilde{\mcl{B}}^*\mcl{P} \quad \text{and}\quad [\mcl{P}\mcl{A}+\mcl{W}\mcl{C}_{q}]=\mcl{P}\tilde{\mcl{A}},
	\end{equation*}
	and therefore
	\begin{align*}
		&Q:=\bmat{-\gamma I & \tilde{\mcl{D}} &~ \tilde{\mcl{C}}\\
			\tilde{\mcl{D}}^* & -\gamma I &~ \tilde{B}^*\mcl{P}\mcl{T}	\\
			\tilde{\mcl{C}}^* & \mcl{T}^*\mcl{P}\tilde{\mcl{B}} &~ \tilde{\mcl{A}}^*\mcl{P}\mcl{T} + \mcl{T}^*\mcl{P}\tilde{\mcl{A}}}
		\preceq 0.
	\end{align*}
	Thus, the conditions of Lem.~\ref{lem:KYP} are satisfied for the PIE~\eqref{eq:PIE_error}, and the result follows. 
\end{pf}

Cor.~\ref{cor:KYP_estimator} proves that, if the LPI~\eqref{eq:KYP_estimator_LPI} is feasible for some $(\gamma,\mcl{P},\mcl{W})$, then, using the estimator defined by~\eqref{eq:PIE_estimator} with gain operator $\mcl{L}:=\mcl{P}^{-1}\mcl{W}$, the $H_{\infty}$-norm $\sup_{w\neq 0}\frac{\|\hat{z}-z\|_{L_2}}{\|w\|_{L_2}}$ of the associated error dynamics is upper-bounded by $\gamma$. Using this result, we finally prove Thm.~\ref{thm:PDE_estimator}

\begin{pf}[Proof of Thm.~\ref{thm:PDE_estimator}]
	Suppose that the conditions of the theorem are satisfied. Fix an arbitrary initial state $\mbf{u}_{0}\in X$ and bounded input $w\in L_{2}^{n_{w}}[0,\infty)$, and let $(\mbf{u},z,\mbf{q})$ be a corresponding solution to the PDE defined by parameters $\tnf{G}_{\tnf{pde}}$. Then, by Lem~\ref{lem:PDE2PIE_observed}, $(\mbf{v},\smallbmat{z\\\mbf{q}})$ with $\mbf{v}=\partial_{x}^2\partial_{y}^2\mbf{u}$ is a solution to the PIE defined by $\{\mcl{T},\mcl{A},\mcl{B},\smallbmat{\mcl{C}\\\mcl{C}_{q}},\smallbmat{\mcl{D}\\\mcl{D}_{q}}\}$, with initial state $\mbf{v}_{0}=\partial_{x}^2\partial_{y}^2\mbf{u}_{0}$. 
	Let $(\hat{\mbf{v}},\hat{z})$ be a solution to the PIE~\eqref{eq:PIE_estimator}
	with input $\mbf{q}(t)$ and initial state $\hat{\mbf{v}}_{0}=\partial_{x}^2\partial_{y}^2\mbf{u}_{0}$, and define $\mbf{e}=\hat{\mbf{v}}-\mbf{v}$ and $\tilde{z}=\hat{z}-z$. Then, 
	$\mbf{e}(0)=\hat{\mbf{v}}_{0}-\mbf{v}_{0}=\mbf{0}$, and $(w(t),\tilde{z}(t))$ satisfies~\eqref{eq:PIE_error} for any $t\geq 0$. By Cor.~\ref{cor:KYP_estimator} we conclude that $\tilde{z}\in L_{2}^{n_{z}}[0,\infty)$ and $\|\tilde{z}\|_{L_2}\leq\gamma \|w\|_{L_2}$.
\end{pf}


\section{Estimator Synthesis of 2D PDEs using Convex Optimization}\label{sec:LPI_implementation}

In the previous section, it was shown that an estimator for a class of 2D PDEs can be synthesized by solving the LPI~\eqref{eq:KYP_estimator_LPI}. Specifically, if for some $\gamma>0$, $\{\mcl{P},\mcl{W}\}$ is a solution to the LPI, then an estimator with $H_{\infty}$-norm of the error bounded by $\gamma$ may be defined as in~\eqref{eq:PIE_estimator}, with gain $\mcl{L}=\mcl{P}^{-1}\mcl{W}$. In this section, we show how this LPI can be numerically solved by parameterizing the PI operator variables $\mcl{P}$ and $\mcl{W}$ by matrices.
In particular, for some $p,r,m\in\N$, fix $\mcl{Z}_{1}\in\PIset_{2D}^{p\times n_{u}}$, $\mcl{Z}_{2}\in\PIset_{2D}^{r\times n_{u}+n_{z}+n_{w}}$, and $\mcl{Z}_{3}\in\PIset_{2D\from 0}^{m\times n_{q_{1}}}\times \PIset_{2D\from 1D}^{m\times n_{q_{2}}+n_{q_{3}}}$ to be defined by monomials of degrees at most $d_{1},d_{2},d_{3}\in\N$, in the variables $x,y,\theta,\eta\in[0,1]$. Then, parameterizing
\begin{equation*}
	\mcl{P}=\mcl{Z}_{1}^* \tnf{M}[P]\mcl{Z}_{1},\quad
	\mcl{Q}=\mcl{Z}_{2}^*\tnf{M}[Q]\mcl{Z}_{2},\enspace\text{and}\enspace
	\mcl{W}=\tnf{M}[W]\mcl{Z}_{3},
\end{equation*}
by matrices $P\in\R^{p\times p}$, $Q\in\R^{r\times r}$, and $W\in\R^{n_{u}\times m}$, it can be shown that $P\succeq 0$ and $Q\preceq 0$ imply $\mcl{P}\succeq 0$ and $\mcl{Q}\preceq 0$, respectively (see e.g.~\cite{jagt2022PIE_2DHinfty}). 
In this manner, the LPI conditions in Thm.~\ref{thm:PDE_estimator} can be enforced as LMI conditions, allowing an $H_{\infty}$-optimal estimator for a 2D PDE to be synthesized as in Algorithm~\ref{alg:estimator}.
\begin{algorithm2e}
	\SetAlgoLined
	\KwData{PDE $\tnf{G}_{\tnf{pde}}$, PI operators $\mcl{Z}_{1},\!\mcl{Z}_{2},\!\mcl{Z}_{3}$, scalar $\epsilon\!>\!0$.}
	1. Compute $\{\mcl{T},\mcl{A},\mcl{B},\mcl{C},\mcl{D},\mcl{C}_{q},\mcl{D}_{q}\}$ as per Thm.~\ref{thm:PDE_estimator}\;
	2. Solve the semidefinite program \\[-1.2em]
	\begin{flalign}\label{eq:estimator_LMI}
		&\mat{\min\limits_{\gamma>0,P,Q,W}\gamma\\
			P\succeq 0,\\[0.3em]Q\preceq 0,}	&
		&\bmat{-\gamma I &\! -\mcl{D} &\! \mcl{C}\\
			(\cdot)^* &\! -\gamma I &\! -[\mcl{B}^*\mcl{P}+\mcl{D}_{q}^*\mcl{W}^*]\mcl{T}	\\
			(\cdot)^* &\! (\cdot)^* &\! (\cdot)^* + \mcl{T}^*[\mcl{P}\mcl{A}+\mcl{W}\mcl{C}_{q}]}
		=\mcl{Q},	& \notag \\[-1.1em]
	\end{flalign}
	where $\mcl{P}=\mcl{Z}_{1}^*\text{M}[P]\mcl{Z}_{1}+\epsilon \tnf{M}[I_{n_{u}}]$, $\mcl{Q}=\mcl{Z}_{2}^*\tnf{M}[Q]\mcl{Z}_{2}$ and $\mcl{W}=\tnf{M}[W]\mcl{Z}_{3}$\;
	3. Compute the Luenberger gain $\mcl{L}=\mcl{P}^{-1}\mcl{W}$\;
	\caption{$H_{\infty}$-Optimal Estimator Synthesis}\label{alg:estimator}
\end{algorithm2e}
\vspace*{-0.2cm} 

Note that, since the operator $\mcl{P}$ in Algorithm~\ref{alg:estimator} is bounded, linear, and coercive, the inverse $\mcl{P}^{-1}$ in $\mcl{L}=\mcl{P}^{-1}\mcl{W}$ is well-defined. In order to actually compute this inverse, we recall that an explicit expression for the operator inverse has already been derived in \cite{miao2019PI_inversion} for a class of \textit{separable} 1D PI operators, taking the form $(\mcl{R}\mbf{v})(x)=R_{0}(x)\mbf{v}(x)+Z(x)^T\int_{0}^{1}HZ(\theta)\mbf{v}(\theta)d\theta$. The following proposition replicates this result for separable 2D PI operators, providing an expression for the inverse of any 2D PI operator $\mcl{R}$ that admits a decomposition $\mcl{R}=\text{M}[R_{0}]+\text{M}[Z]H\smallint_{[0,1]^2}[Z]$.
	
\begin{prop}\label{prop:PI_inverse}
	For given $n\in\N$, let $R_{0}\in L_2^{n\times n}[[0,1]^2]$ be invertible, and let $R_{1}\in L_2^{n\times n}[[0,1]^2\times[0,1]^2]$ be such that
	\begin{align*}
		R_{1}(x,y,\theta,\eta)\!&=\!Z(x,y)^T HZ(\theta,\eta),	&	x,y,\theta,\eta&\!\in[0,1],
	\end{align*}
	for some $Z\in L_2^{p\times n}[[0,1]^2]$ and $H\in\R^{p\times p}$ with $p\in\N$. Let
	\begin{align*}
		Q_{0}(x,y)&:=[R_{0}(x,y)]^{-1},	\\
		Q_{1}(x,y,\theta,\eta)&:=Q_0(x,y)Z(x,y)^T\hat{H}Z(\theta,\eta)Q_0(\theta,\eta),
	\end{align*}
	for $x,y,\theta,\eta\in[0,1]$, where 
	\begin{equation*}
		\hat{H}:=-H(I_{p}+KH)^{-1} =-(I_{p}+HK)^{-1}H\in\R^{p\times p}
	\end{equation*}
	with $K:=\int_{0}^{1}\int_{0}^{1}Z(\nu,\mu)Q_0(\nu,\mu)Z(\nu,\mu)^Td\mu d\nu\in\R^{p\times p}$. If $R:=\smallbmat{R_{0}&0&0\\0&R_{1}&R_{1}\\0&R_{1}&R_{1}}$ and $Q:=\smallbmat{Q_{0}&0&0\\0&Q_{1}&Q_{1}\\0&Q_{1}&Q_{1}}$, then
	\begin{equation*}
		\PI[Q]\circ\PI[R]=\PI[R]\circ\PI[Q]=\text{M}[I_{n}].
	\end{equation*}
\end{prop}

\begin{pf}
	Fix arbitrary $R_{0}\in L_2^{n\times n}[[0,1]^2]$ and $R_{1}\in L_2^{n\times n}[[0,1]^2\times[0,1]^2]$ as proposed, and let associated $Q_{0}$ and $Q_{1}$ be as defined. Then, by the composition rules of 2D PI operators, we note that
	\begin{align*}
		\PI[Q]\circ\PI[R]=\PI\smallbmat{Q_{0}&0&0\\0&Q_{1}&Q_{1}\\0&Q_{1}&Q_{1}}\circ\PI\smallbmat{R_{0}&0&0\\0&R_{1}&R_{1}\\0&R_{1}&R_{1}}=\PI\smallbmat{P_{0}&0&0\\0&P_{1}&P_{1}\\0&P_{1}&P_{1}},
	\end{align*}	
	where
	\begin{align*}
		&P_{0}(x,y)\!=Q_{0}(x,y)R_{0}(x,y)=I_{n},	\\
		&P_{1}(x,y,\theta,\eta)\!=\!Q_{0}(x,y)R_{1}(x,y,\theta,\eta)\!+\!Q_{1}(x,y,\theta,\eta)R_{0}(\theta,\eta) \\
		&\quad\hspace*{1.0cm} + \int_{0}^{1}\int_{0}^{1}Q_1(x,y,\nu,\mu)R_1(\nu,\mu,\theta,\eta)d\mu d\nu	\\
		&\quad=Q_{0}(x,y)Z(x,y)^T HZ(\theta,\eta)\\
		&\quad\hspace*{0.5cm}+Q_{0}(x,y)Z(x,y)^T\hat{H}Z(\theta,\eta)Q_{0}(\theta,\eta)R_{0}(\theta,\eta) 	\\
		&\quad\hspace*{1.0cm}+Q_{0}(x,y)Z(x,y)^T\hat{H}K HZ(\theta,\eta)	\\
		&\quad=Q_{0}(x,y)Z(x,y)^T\bl[H+\hat{H}+\hat{H}KH\br]Z(\theta,\eta)	\\
		&\quad=Q_{0}(x,y)Z(x,y)^T\bl[\hat{H}(I_{p}+KH)+H\br]Z(\theta,\eta)=0.
	\end{align*}
	It follows that, $\PI[Q]\circ\PI[R]=\PI\smallbmat{I_{n}&0&0\\0&0&0\\0&0&0}=\text{M}[I_{n}]$. Similarly, we find that
	\begin{align*}
		\PI[R]\circ\PI[Q]=\PI\smallbmat{R_{0}&0&0\\0&R_{1}&R_{1}\\0&R_{1}&R_{1}}\circ \PI\smallbmat{Q_{0}&0&0\\0&Q_{1}&Q_{1}\\0&Q_{1}&Q_{1}}=\PI\smallbmat{P_{0}&0&0\\0&P_{1}&P_{1}\\0&P_{1}&P_{1}},
	\end{align*}
	where
	\begin{align*}
		&P_{0}(x,y)\!=R_{0}(x,y)Q_{0}(x,y)=I_{n},	\\
		&P_{1}(x,y,\theta,\eta)\!=\!R_{0}(x,y)Q_{1}(x,y,\theta,\eta)\!+\!R_{1}(x,y,\theta,\eta)Q_{0}(\theta,\eta) \\
		&\quad\hspace*{1.0cm} + \int_{0}^{1}\int_{0}^{1}R_1(x,y,\nu,\mu)Q_1(\nu,\mu,\theta,\eta)d\mu d\nu	\\
		&\quad=R_{0}(x,y)Q_{0}(x,y)Z(x,y)^T \hat{H}Z(\theta,\eta)Q_{0}(\theta,\eta)\\
		&\quad\hspace*{0.5cm}+Z(x,y)^T HZ(\theta,\eta)Q_{0}(\theta,\eta) 	\\
		&\quad\hspace*{1.0cm}+Z(x,y)^TH K \hat{H}Z(\theta,\eta)Q_{0}(\theta,\eta)	\\
		&\quad=Z(x,y)^T\bl[\hat{H}+H+HK\hat{H}\br]Z(\theta,\eta)Q_{0}(\theta,\eta)	\\
		&\quad=Z(x,y)^T\bl[(I_{p}+HK)\hat{H}+H\br]Z(\theta,\eta)Q_{0}(\theta,\eta)=0.
	\end{align*}
	Thus, also $\PI[R]\circ\PI[Q]=\text{M}[I_{n}]$, concluding the proof.	
%
\end{pf}

Using Prop.~\ref{prop:PI_inverse}, the operator $\mcl{P}=\mcl{Z}_{1}^*\tnf{M}[P]\mcl{Z}_{1}+\epsilon\tnf{M}[I_{n_{u}}]$ can be inverted analytically if it is chosen to have a suitable -- separable -- structure. Here, since any polynomial function $R\in L_{2}[[0,1]^4]$ is separable as $R(x,y,\theta,\eta)=Z(x,y)^T HZ(\eta,\theta)$ for some matrix $H$ and vector of monomials $Z$, the operator $\mcl{P}$ can be inverted by Prop.~\ref{prop:PI_inverse} if it is chosen to be of the form
\begin{equation*}
	(\mcl{P}\mbf{v})(x,y)\!=\!R_{0}(x,y)\mbf{v}(x,y) +\!\int_{0}^{1}\!\!\int_{0}^{1}\!\!R_{1}(x,y,\theta,\eta)\mbf{v}(\theta,\eta)d\theta d\eta.
\end{equation*}
for some polynomials $R_{0},R_{1}$. 
Of course, by imposing such a restriction on the structure of $\mcl{P}$, as well as allowing $\mcl{P}$ and $\mcl{W}$ to be defined only by polynomial functions in general, we are necessarily introducing conservatism.
However, the fact that both $\mcl{P}$ and $\mcl{W}$ in~\eqref{eq:estimator_LMI} may be defined by (partial) integral operators still allows significantly more freedom than parameterizing the Luenberger gain by merely a multiplier operator, as is commonly done in practice. Moreover, conservatism can be reduced by increasing the maximal degrees of the monomials defining each operator variable.

\begin{figure*}[t!]
	\begin{subfigure}[b]{0.516\linewidth}
		\centering
		\hspace*{-0cm}\includegraphics[width=0.99\linewidth]{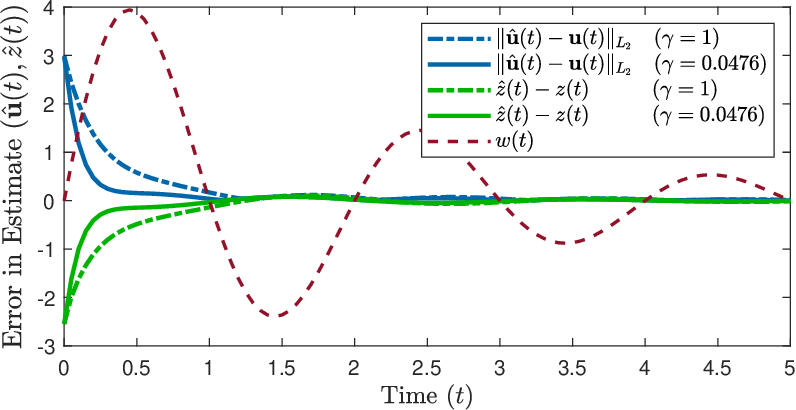}
		\caption{\footnotesize 
			$r=4$}
		\label{fig:Ex1:r4}
	\end{subfigure}
	\begin{subfigure}[b]{0.479\linewidth}
		\centering
		\hspace*{-0cm}\includegraphics[width=0.99\linewidth]{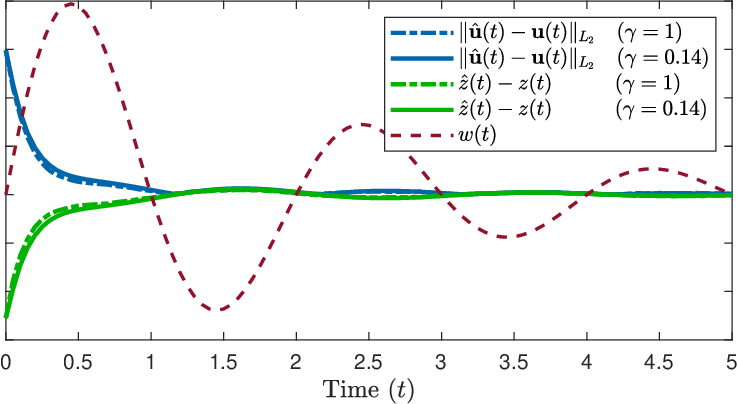}
		\caption{\footnotesize 
			$r=8$}
		\label{fig:Ex1:r8}
	\end{subfigure}
	\caption{\footnotesize 
		Error in the estimate of the PDE state $\mbf{u}(t)$ and the output $z(t)$ for the PDE in Subsection~\ref{sec:Numerical_Examples:1} with $r=4$ and $r=8$, for disturbance $w(t)=5e^{-t/2}\sin(\pi t)$ and initial error $\hat{\mbf{u}}(0,x,y)-\mbf{u}(0,x,y)=5((x-1)^4-1)\sin(0.5\pi y)$. Errors are plotted for estimators achieved both for fixed $\gamma=1$ in the LMI~\eqref{eq:estimator_LMI}, as well as minimizing over the value of $\gamma$ in this LMI.
	}
	\label{fig:Ex1}
\end{figure*}

\section{A Numerical Example}\label{sec:Numerical_Examples}

The presented methodology for estimator synthesis of 2D PDEs has been fully incorporated into the PIETOOLS software suite (see~\cite{shivakumar2021PIETOOLS}). This software offers an interface for declaring a desired 1D or 2D PDE as a structure \texttt{PDE}, at which point it can be converted to an associated PIE by calling \texttt{PIE=convert(PDE)}. Then, the LMI~\eqref{eq:estimator_LMI} can be constructed and solved by calling \texttt{[Lop,gam]=lpisolve(PIE,`estimator')}, returning an object \texttt{Lop} representing the optimal gain $\mcl{L}=\mcl{P}^{-1}\mcl{W}$, and a scalar \texttt{gam} representing this lowest upper bound $\gamma$ on the $H_{\infty}$-norm of the associated error dynamics.

In this section, we use the PIETOOLS software to construct an optimal estimator for the unstable 2D heat equation in~\eqref{eq:PDE_example_intro}, observing the PDE state along the upper boundary of the domain. To test performance of the obtained estimator, we simulate the error $\mbf{e}(t)$ in the fundamental state $\mbf{v}(t)$ based on the PIE~\eqref{eq:PIE_error}, using a Galerkin method with 36 basis functions for the spatial discretization, and an explicit Euler scheme with time step $\Delta t=2\cdot 10^{-4}$ for the temporal integration. 
More details on this simulation scheme can be found in Appendix~\ref{appx:Simulation}.

\subsection{Estimator for a 2D Heat Equation}\label{sec:Numerical_Examples:1}

We consider the following 2D heat equation, with sensors along the upper boundary,
\begin{align*}
	\mbf{u}_{t}(t)&=\mbf{u}_{xx}(t)+\mbf{u}_{yy}(t)+r\mbf{u}(t) +w(t),	\\
	z(t)&=\int_{0}^{1}\int_{0}^{1}\mbf{u}(t,x,y)dx dy,	\notag\\
	\mbf{q}_{1}(t)&=\mbf{u}(t,1,\cdot)+\eta_{1}(t),\hspace*{0.45cm} \mbf{q}_{2}(t)=\mbf{u}(t,\cdot,1)+\eta_{2}(t),	\notag\\
	\mbf{u}(t,0,\cdot)&=\mbf{u}_{x}(t,1,\cdot)\equiv 0,\quad \mbf{u}(t,\cdot,0)=\mbf{u}_{y}(t,\cdot,1)\equiv 0,
\end{align*}
where $\mbf{u}(t)\in L_2[[0,1]^2]$, $w(t),z(t),\eta_{1}(t),\eta_{2}(t)\in\R$, and $\mbf{q}(t)=\smallbmat{\mbf{q}_{1}(t)\\\mbf{q}_{2}(t)}\in L_2^2[0,1]$. Using PIETOOLS, we construct a PIE representation associated to this system and synthesize an estimator for values of the parameter $r=4$ and $r=8$, for which the system is stable and unstable, respectively. Solving the optimization program in~\eqref{eq:estimator_LMI}, we obtain optimal gain operators $\mcl{L}$ that achieve lowest bounds on the $H_{\infty}$-norm of the associated error dynamics as $\gamma=0.0476$ and $\gamma=0.1403$ for $r=4$ and $r=8$, respectively. In addition, for each value of $r\in\{4,8\}$, sub-optimal gains $\mcl{L}$ were also computed for a fixed bound $\gamma=1$ in the LMI.

The error dynamics corresponding to the obtained estimators were simulated with disturbance $w(t)=5e^{-t/2}\sin(\pi t)$ and initial state $\mbf{u}_{0}(t,x,y)=5((x-1)^4-1)\sin(0.5\pi y)$, starting with an initial state estimate $\hat{\mbf{u}}(0)=0$. Sensor noise $\eta_{1}(t)$ and $\eta_{2}(t)$ at each time was generated from a Gaussian distribution with mean $0$ and variance $0.04$. Figure~\ref{fig:Ex1} shows the norm of the error in the PDE state and output for $t\in[0,5]$, as well as the value of the disturbance. The plots show that, both for the stable ($r=4$) and unstable ($r=8$) PDE, the error in the estimates of the state and output converge to zero, with similar or faster convergence achieved using the optimal estimators associated to the minimized bounds $\gamma$.



	\section{Conclusion}
	
	In this paper, a new convex-optimization-based method was presented for estimator synthesis of linear, 2nd order, 2D PDEs with state observations along the boundary. To this end, it was proved that any sufficiently well-posed such PDE can be equivalently represented as a PIE, specifically proving that the value of the state $\mbf{u}$ along the boundary can be expressed in terms of a PI operator acting on the fundamental state $\partial_{x}^2\partial_{y}^2\mbf{u}$. Parameterizing a Luenberger-type estimator for the PIE by a PI operator, it was then shown that a value of this operator with guaranteed bound on the $H_{\infty}$-norm of the estimator error dynamics can be computed by solving an LPI, which in turn could be solved as an LMI. The proposed methodology has been incorporated in the PIETOOLS software suite, and applied to construct an estimator for an unstable heat equation, using simulation to show convergence of the estimated state to the true value.
	

\fontsize{10pt}{5pt}\selectfont
\bibliography{bibfile}

\begin{thebibliography}{15}
\providecommand{\natexlab}[1]{#1}
\providecommand{\url}[1]{\texttt{#1}}
\providecommand{\urlprefix}{URL }
\expandafter\ifx\csname urlstyle\endcsname\relax
  \providecommand{\doi}[1]{doi:\discretionary{}{}{}#1}\else
  \providecommand{\doi}{doi:\discretionary{}{}{}\begingroup
  \urlstyle{rm}\Url}\fi

\bibitem[{Das et~al.(2019)Das, Shivakumar, Weiland, and
  Peet}]{das2019PIE_estimation}
Das, A., Shivakumar, S., Weiland, S., and Peet, M.M. (2019).
\newblock {$\mcl{H}_{\infty}$} optimal estimation for linear coupled {PDE}
  systems.
\newblock In \emph{2019 IEEE 58th Conference on Decision and Control (CDC)},
  262--267. IEEE.

\bibitem[{Duan and Yu(2013)}]{duan2013LMIs}
Duan, G.R. and Yu, H.H. (2013).
\newblock \emph{{LMI}s in control systems: analysis, design and applications}.
\newblock CRC press.

\bibitem[{Hasan et~al.(2016)Hasan, Aamo, and
  Krstic}]{hasan2016BacksteppingObserver_PDEODE}
Hasan, A., Aamo, O.M., and Krstic, M. (2016).
\newblock Boundary observer design for hyperbolic {PDE--ODE} cascade systems.
\newblock \emph{Automatica}, 68, 75--86.

\bibitem[{Holmes et~al.(1994)Holmes, Lewis, Banks, and
  Veit}]{holmes1994partial}
Holmes, E.E., Lewis, M., Banks, J., and Veit, R.R. (1994).
\newblock Partial differential equations in ecology: Spatial interactions and
  population dynamics.
\newblock \emph{Ecology}, 75, 17--29.

\bibitem[{Jadachowski et~al.(2015)Jadachowski, Meurer, and
  Kugi}]{jadachowski2015BacksteppingObserver_ND}
Jadachowski, L., Meurer, T., and Kugi, A. (2015).
\newblock Backstepping observers for linear {PDE}s on higher-dimensional
  spatial domains.
\newblock \emph{Automatica}, 51, 85--97.

\bibitem[{Jagt and Peet(2021)}]{jagt2021PIEArxiv}
Jagt, D.S. and Peet, M.M. (2021).
\newblock A {PIE} representation of coupled {2D} {PDE}s and stability analysis
  using {LPI}s.
\newblock \emph{arXiv eprint:2109.06423}.

\bibitem[{Jagt and Peet(2022)}]{jagt2022PIE_2DHinfty}
Jagt, D.S. and Peet, M.M. (2022).
\newblock {$L_2$}-gain analysis of coupled linear {2D} {PDE}s using linear {PI}
  inequalities.
\newblock In \emph{2022 IEEE 61st Conference on Decision and Control (CDC)},
  6097--6104.

\bibitem[{Lhachemi and Prieur(2022)}]{lhachemi2021outputfeedackstabilization}
Lhachemi, H. and Prieur, C. (2022).
\newblock Boundary output feedback stabilization of a class of
  reaction-diffusion {PDE}s with delayed boundary measurement.
\newblock \emph{International Journal of Control}, 1--11.

\bibitem[{Meurer(2013)}]{meurer2013BacksteppingObserver_Semilinear}
Meurer, T. (2013).
\newblock On the extended {Luenberger}-type observer for semilinear
  distributed-parameter systems.
\newblock \emph{IEEE Transactions on Automatic Control}, 58(7), 1732--1743.

\bibitem[{Miao et~al.(2019)Miao, Peet, and Gu}]{miao2019PI_inversion}
Miao, G., Peet, M.M., and Gu, K. (2019).
\newblock Inversion of separable kernel operator and its application in control
  synthesis.
\newblock \emph{Delays and Interconnections: Methodology, Algorithms and
  Applications}, 265--280.

\bibitem[{Shivakumar et~al.(2021)Shivakumar, Das, and
  Peet}]{shivakumar2021PIETOOLS}
Shivakumar, S., Das, A., and Peet, M. (2021).
\newblock {PIETOOLS} 2020a: User manual.
\newblock \emph{arXiv preprint arXiv:2101.02050}.

\bibitem[{Wang and Fridman(2024)}]{wang2024Observer_Stochastic}
Wang, P. and Fridman, E. (2024).
\newblock Sampled-data finite-dimensional observer-based control of {1D}
  stochastic parabolic {PDE}s.
\newblock \emph{SIAM Journal on Control and Optimization}, 62(1), 297--325.

\bibitem[{Yu et~al.(2020)Yu, Gan, Bayen, and
  Krstic}]{yu2020BacksteppingObserver_traffic}
Yu, H., Gan, Q., Bayen, A., and Krstic, M. (2020).
\newblock {PDE} traffic observer validated on freeway data.
\newblock \emph{IEEE Transactions on Control Systems Technology}, 29(3),
  1048--1060.

\bibitem[{Zayats et~al.(2021)Zayats, Fridman, and Zhuk}]{zayats2021ObserverNS}
Zayats, M., Fridman, E., and Zhuk, S. (2021).
\newblock Global observer design for {Navier-Stokes} equations in {2D}.
\newblock In \emph{2021 60th IEEE Conference on Decision and Control (CDC)},
  1862--1867.

\bibitem[{Zuazua(2005)}]{zuazua2005Observability_FiniteDifference}
Zuazua, E. (2005).
\newblock Propagation, observation, and control of waves approximated by finite
  difference methods.
\newblock \emph{SIAM review}, 47(2), 197--243.

\end{thebibliography}
	
\bigskip
	
\appendix
\section{Numerical Simulation of PIEs}\label{appx:Simulation}

In section~\ref{sec:Numerical_Examples}, we applied the proposed LPi methodology to construct an estimator for a 2D heat equation with state observations along the boundary. To test the performance of this estimator, we simulated the dynamics of the associated error $\mbf{e}(t)=\mbf{v}(t)-\hat{\mbf{v}}(t)$ in the fundamental state $\mbf{v}(t)=\partial_{x}^2\partial_{y}^2\mbf{u}(t)$, defined by the PIE~\eqref{eq:PIE_estimator}. In particular, we model this error using a Galerkin approach, wherein we consider the error state $\mbf{e}(t)$ only on a finite-dimensional subspace $Y\subseteq L_2[[0,1]]$, spanned by a basis of pairwise orthonormal functions $\mbs{\phi}_{i}\in L_2[[0,1]^2]$ for $i\in\{1,\hdots,N\}$, so that $\ip{\mbs{\phi}_{i}}{\mbs{\phi}_{j}}_{L_2}=1$ if $i=j$ but $0$ otherwise. Then, the error at any time is defined by a vector of coefficients $c_{k}(t)\in\R^{N}$ as $\mbf{e}(t)=\sum_{j=1}^{N}c_{j}(t)\mbs{\phi}_{j}$, and we can express
\begin{equation*}
	\mcl{T}\mbf{e}_{t}(t)=\sum_{j=1}^{N}\dot{c}_{j}(t) (\mcl{T}\mbs{\phi}_{j}),
	\qquad
	\bar{\mcl{A}}\mbf{e}(t)=\sum_{j=1}^{N}c_{j}(t) (\bar{\mcl{A}}\mbs{\phi}_{j}).
\end{equation*}
We require this error to weakly satisfy the PIE~\eqref{eq:PIE_error} on the space $Y$, so that
\begin{equation*}
	\ip{\mbs{\phi}_{i}}{\mcl{T}\mbf{e}_{t}(t)}_{L_2}=\ip{\mbs{\phi}_{i}}{\bar{\mcl{A}}\mbf{e}(t) +\bar{\mcl{B}}w(t)},\quad \forall i\in\{1,\hdots,N\}.
\end{equation*}
Defining matrices $T,\bar{A}\in\R^{N\times N}$ by $T_{ij}:=\ip{\mbs{\phi}_{i}}{\mcl{T}\mbs{\phi}_{j}}_{L_2}$ and $\bar{A}_{ij}:=\ip{\mbs{\phi}_{i}}{\bar{\mcl{A}}\mbs{\phi}_{j}}_{L_2}$, the resulting system can then be represented as an ODE on the vector of coefficients $c(t)\in\R^{N}$
\begin{equation*}
	T\dot{c}(t)=Ac(t) +\ip{\mbs{\phi}}{\mcl{B}w(t)}_{L_2}.
\end{equation*}
Assuming the matrix $T$ to be invertible, this system can be numerically solved for a given input $w$ using any time-stepping scheme, starting with an initial value of the coefficients as $c_{j}(0)=\ip{\mbs{\phi}_{j}}{\mbf{e}(0)}_{L_2}=\ip{\mbs{\phi}_{j}}{\partial_{x}^2\partial_{y}^2\hat{\mbf{u}}(0)-\partial_{x}^2\partial_{y}^2\mbf{u}(0)}_{L_2}$. 
From the obtained values of the coefficients at each time, the error in the estimate of the fundamental state can then be computed as $\mbf{e}(t)=\sum_{j=1}^{N}c_{j}(t)\mbs{\phi}_{j}$, and we can retrieve the error in the estimate of the PDE state as $\mcl{T}\mbf{e}(t)=\sum_{j=1}^{N}c_{j}(t)\mcl{T}\mbs{\phi}_{j}$. The error $\hat{z}(t)$ in the regulated output can then be computed directly from the PIE~\eqref{eq:PIE_error}.

For the numerical implementation in Section~\ref{sec:Numerical_Examples}, the basis functions $\mbs{\phi}_{i}$ were chosen as Legendre polynomials of degree at most $5$ in each variable $x,y\in[0,1]$, yielding a total of $N=(5+1)^2=36$ basis functions. For the numerical time-stepping, an explicit Euler scheme was applied as
\begin{equation*}
	c(t_{k+1})=c(t_{k}) +\Delta t(T^{-1}Ac(t_{k}) +\ip{\mbs{\phi}}{\mcl{B}w(t)}_{L_2})
\end{equation*}
using a time step $\Delta t=2\cdot 10^{-4}$, and simulating from $t_{0}=0$ up to $t_{N_{t}}=10$.
	
\end{document}